\documentclass[11pt,letterpaper]{amsart}
\usepackage{graphicx}

\usepackage{amsmath}
\usepackage{amsfonts}
\usepackage{enumerate}
\usepackage{amssymb}
\usepackage{epstopdf}
\usepackage{tikz-cd}

\usepackage{anysize}
\marginsize{3cm}{3cm}{3cm}{3cm}

\newtheorem{theorem}{Theorem}[section]
\newtheorem{proposition}[theorem]{Proposition}
\newtheorem{lemma}[theorem]{Lemma}

\newtheorem{conjecture}[theorem]{Conjecture}

\theoremstyle{definition}

\newtheorem*{definition*}{Definition}

\numberwithin{equation}{section}

\renewcommand{\epsilon}{\varepsilon}

\renewcommand{\rho}{\varrho}

\long\def\forget#1\forgotten{}
\begin{document}

\title{Small separators, upper bounds for $l^\infty$-widths,  and systolic geometry}

\author{Sergey Avvakumov, Alexander Nabutovsky}
\date{}
\maketitle

\begin{abstract} 


The goal of this paper is to prove the first estimates in systolic geometry with a sublinear dependence on the dimension.
We prove that if $M^n\subset \mathbb{R}^N$ is {\it essential}, then there exists a non-contractible closed curve on $M^n$ contained in a cube in $\mathbb{R}^N$ with side length $const\ \sqrt{n}\ vol^{\frac{1}{n}}(M^n)$ with sides parallel to the coordinate axes. If the codimension is $1$, then the side length of the cube is $4\ vol^{\frac{1}{n}}(M^n)$. 

This result follows from the following inequalities that relate the Euclidean volume of a closed submanifold $M^n\subset \mathbb{R}^N$ with its width $W^{l^\infty}_{n-1}(M^n)$ with respect to the $l^\infty$-metric on the ambient $\mathbb{R}^N$. (This width is defined as the infimum over all continuous maps $\phi:M^n\longrightarrow K^{n-1}\subset\mathbb{R}^N$ of $sup_{x\in M^n}\Vert \phi(x)-x\Vert_{l^\infty}$.)
We prove that $W^{l^\infty}_{n-1}(M^n)\leq const\ \sqrt{n}\ vol(M^n)^{\frac{1}{n}}$, and if the codimension $N-n$ is equal to $1$, then $W^{l^\infty}_{n-1}(M^n)\leq \sqrt{3}\ vol(M^n)^{\frac{1}{n}}$.

Our proof of these results uses several new ideas. First, we prove a version of the kinematic formula, where one averages over isometries of $l^N_\infty$ (Theorem 3.5). Then
we construct high-codimension
analogs of optimal foams discovered in [KORW] and [AK] (Theorem 3.9). Finally, we develop a new approach in systolic geometry that can be described as a non-linear version of the classical Federer-Fleming argument, where we construct continuous maps of $M^n\subset \mathbb{R}^N$ to lower-dimensional polyhedra by pushing out from 
specially designed non-linear $(N-n)$-dimensional complexes in $\mathbb{R}^N$ that do not intersect $M^n$.

\end{abstract}

\section{Introduction}

\subsection{Main Results: new width-volume inequalities}

Let $M^n\subset \mathbb{R}^N$ be a closed manifold.  Denote by $vol(M^n)$ the Euclidean volume of $M^n$.

\begin{definition*}
The $l^\infty$-width $W_{n-1}^{l^\infty}(M^n)$ of $M^n$ is the infimum over all $(n-1)$-dimensional polyhedra $K^{n-1}\subset\mathbb{R}^N$ and all continuous maps $\phi:M^n\longrightarrow K^{n-1}$
of $\sup_{x\in M^n}\Vert x-\phi(x)\Vert_{l^\infty}$, where the distance between $x$ and $\phi(x)$ is calculated with respect to the $l^\infty$-norm.
\end{definition*}

Informally,  $W_{n-1}^{l^\infty}(M^n)$ measures how far $M^n$ is from an $(n-1)$-dimensional polyhedron. 

Here are our main results:
\begin{theorem}
Let $M^n\subset \mathbb{R}^N$ be a closed manifold.
Then $$W^{l^\infty}_{n-1}(M^n)\leq C\ \sqrt{n}\ vol(M^n)^{\frac{1}{n}},$$
where $C$ is an absolute constant (which, in particular, does not depend on $n$). We can take $C=\frac{2\pi}{\sqrt{e}}(1+o(1))$ as $n\longrightarrow\infty$.
\end{theorem}

\begin{theorem}
Let $M^n\subset \mathbb{R}^{n+1}$ be a closed hypersurface.
Then $$W^{l^\infty}_{n-1}(M^n)\leq 3^{\frac{1}{n}}\cdot vol(M^n)^{\frac{1}{n}}\leq \sqrt{3}\cdot vol(M^n)^{\frac{1}{n}}.$$
\end{theorem}

An example of a smoothed-out boundary of the unit cube demonstrates that Theorem 1.2 is optimal up
to a constant factor. (In other words, one cannot replace the factor $3^{\frac{1}{n}}$ on the right-hand  side of the inequality by $o(1)$, as $n\longrightarrow\infty$).
On the other hand, we do not know if the estimate in Theorem 1.1 is optimal up to a constant factor. 

\subsection{$l^\infty$-width and Urysohn width.}
There are several ways to define the $(n-1)$th ``width'' of a manifold (or more generally a space). Each definition measures how far the manifold is from an $(n-1)$-dimensional polyhedron and they differ only in the way ``far'' is defined.

\begin{definition*}
To define the $(n-1)$-dimensional width of $M^n$ we take the infimum over all $(n-1)$-dimensional polyhedra $K^{n-1}$ and all continuous maps $\phi:M^n\longrightarrow K^{n-1}$ of the following quantity:
\begin{itemize}
    \item For Urysohn width $UW_{n-1}(M^n)$ the quantity is $\sup_{k\in K^{n-1}}diam(\phi^{-1}(k))$.
    \item For $W_{n-1}^X(M^n)$ the quantity is $\sup_{x\in M^n}dist_X(x,\phi(x))$. Here both $M^n$ and $K^{n-1}$ are assumed to be embedded into a Banach space, or a Riemannian manifold, or, more generally, a metric space $X$. One can analogously define $W_{m}^X(M^n)$ by replacing $(n-1)$ with any $m<n-1$.
    \item We write $W^{Eucl}_{m}(M^n)$ for $W_{n-1}^X(M^n)$ in the case when $X$  is $\mathbb{R^N}$ with the standard Euclidean metric.
\end{itemize}
\end{definition*}

Observe that if $k\in K^{n-1}\subset \mathbb{R}^N$, and there exists $\phi: M^n\longrightarrow K^{n-1}$ such that
$\sup_{x\in M^n}\Vert\phi(x)-x\Vert_{l^\infty}\leq W^{l^\infty}_{n-1}(M^n)+\epsilon$,
then for every two points $x_1, x_2\in \phi^{-1}(k)\subset M^n\subset \mathbb{R}^N$ the $l^\infty$-distance between $x_1,x_2$ does not exceed $2W^{l^\infty}_{n-1}(M^n)+2\epsilon$. Therefore, if one considers $M^n$ with the extrinsic metric induced by the ambient $N$-dimensional space regarded as $l^N_\infty$, then we see that $UW_{n-1}(M^n)\leq 2W^{l^\infty}_{n-1}(M^n)$.

If $X$ is the ambient $\mathbb{R}^N$ but regarded as $l^N_\infty$, the width $W_{n-1}^X(M^n)$ becomes $W^{l^\infty}_{n-1}(M^n)$ that was already defined above. 
On the other hand, consider the (distance-preserving) Kuratowski embedding
$M^n\longrightarrow L^\infty(M^n)$ that sends each point $x\in M^n$ to the function $dist_x(y)$ on $M^n$ defined by the formula $dist_x(y)=dist_{M^n}(x,y)$. 
Proposition D in Appendix 1 of [G] implies that the corresponding width $W^{L^\infty(M^n)}(M^n)$ is equal
to $2UW_{n-1}(M^n)$.

Widths $W_m^X(M^n)$ were introduced in Appendix 1 of [G], where they were denoted $Rad_m(M^n\subset X)$.
We refer the reader to Appendix 1 of [G] for
the other information about these metric invariants including their other equivalent definitions. 


\subsection{Essential manifolds and systolic geometry.}

\begin{definition*}
A non-simply connected manifold $M^n$ is called {\it essential} if there is {\it no} continuous map from $M^n$ to a polyhedron $K$ of dimension at most $n-1$ that induces an isomorphism of fundamental groups of $M^n$ and $K$.
\end{definition*}

\begin{definition*}
The \emph{systole}, denoted as $sys_1(M^n)$, of a Riemannian manifold $M^n$ is the length of the shortest non-contractible closed curve on $M^n$. 
\end{definition*}

 Examples of essential closed manifolds include all non-simply connected closed surfaces and, more generally, all aspherical manifolds, as well as all projective spaces $\mathbb{R}P^n$. In his foundational paper, [G]
Gromov proved that $sys_1(M^n)$
satisfies the inequality
$$sys_1(M^n)\leq c(n)\ vol^{\frac{1}{n}}(M^n).$$ This inequality is now known as Gromov's systolic inequality. (On the other hand, I. Babenko ([B]) demonstrated that this inequality does not hold for any non-essential manifold $M^n$.) The original proof of Gromov yielded $c(n)\sim n^{{\frac{3}{2}}n}$. A later, simpler, proof by S. Wenger ([W]) yields $c(n)\sim n^n$. The best currently known upper bound for $c(n)$ is $2({\frac{n!}{2}})^{\frac{1}{n}}\leq n$. This estimate was found in [N] using the ideas from the paper [P] by P. Papasoglu. Papasoglu was not interested in estimating $c(n)$, but construction in [P] (as is) would yield an upper bound of the form
$\exp(const\ n)$. 

On the other hand, the round metrics on $\mathbb{R}P^n$ (as well as some flat metrics
on tori) have $sys_1 ~ \sim\sqrt{n}\ vol^{\frac{1}{n}}$, and Gromov conjectured that the optimal value of $c(n)$ is attained on round projective spaces. We can also consider weaker conjecture that the optimal value of $c(n)$ exceeds $sys_1((\mathbb{R}P^n, can))/vol^{\frac{1}{n}}((\mathbb{R}P^n, can))$ by not more than a constant factor that does not depend on the dimension:

\begin{conjecture}
There exists an absolute constant $C$ such that for each essential
closed Riemannian manifold $M^n$ $sys_1(M^n)\leq C\ \sqrt{n}\ vol^{\frac{1}{n}}(M^n)$.
\end{conjecture}


\subsection{Widths and systole.}
In [G] Gromov conjectured that $$UW_{n-1}(M^n)\leq c(n)\ vol^{\frac{1}{n}}(M^n).\hskip 6truecm (1)$$ This conjecture was proven by L. Guth in [Gu1], [Gu2].
The statement is stronger than the systolic inequality due to the inequality $sys_1(M^n)\leq 3 UW_{n-1}(M^n)$ that holds for all essential manifolds $M^n$. This inequality was proven by Gromov in the same paper [G]. 
The proof is short and not difficult. 
A different short proof of a similar inequality
$sys_1(M^n)\leq 2\ UW_{n-1}(M^n)$ for essential $M^n$
based on the ideas of A. Schwarz ([S]) and observation of R. Karasev that they can be applied in this situation can be found in [N].

It is now natural to conjecture that

\begin{conjecture}
There exists a constant $C$ such that for each Riemannian manifold (and even Riemannian polyhedron) $M^n$
$$UW_{n-1}(M^n)\leq C\ \sqrt{n}\ vol^{\frac{1}{n}}(M^n).$$
\end{conjecture}

Note that Gromov's proof of the inequality $sys_1(M^n)\leq 3UW_{n-1}(M^n)$ for essential $M^n$ as well as the proof of a similar inequality $sys_1(M^n)\leq 2\ UW_{n-1}(M^n)$ in [N] both hold for Riemannian polyhedra.
Therefore, Conjecture 1.4 immediately implies Conjecture 1.3
for Riemannian polyhedra. In [N] it was proven that $$UW_{n-1}(M^n)\leq 2(n!/2)^{\frac{1}{n}}vol^{\frac{1}{n}}(M^n)\leq n\ vol^{\frac{1}{n}}(M^n).$$


One of the central ideas of [G] is to use the Kuratowski embedding $\kappa:M^n\longrightarrow L^{\infty}(M^n)$ that sends each point $x\in M^n$
to the distance function $dist(x, *)$. As this embedding is distance-preserving, the Hausdorff measures
of the images of subsets of $M^n$ in $L^\infty(M^n)$ coincide
with the Hausdorff measures of the same subsets in $M^n$.

It was observed in [Gu0] that there is no need to use
an infinite-dimensional space $L^\infty(M^n)$ here. One can use only the distance functions
to points from a sufficiently dense (finite) set. Although the embedding will not be isometric, the map from $M^n$ to its image in a finite-dimensional $l^N_\infty$ can be made $(1+\epsilon)$-bilipschitz for an arbitrarily small $\epsilon$.

Therefore, Conjectures 1.3 and 1.4 would follow from the positive answer to the following problem: 

\par\noindent
{\bf Problem 1.} {\it Let $M^n\subset l^N_\infty$ be a closed manifold (or a polyhedron). Is it true that
$$W^{l^\infty}_{n-1}(M^n)\leq const\ \sqrt{n}\ H^{\frac{1}{n}}_n(M^n),$$
where $H_n$ denotes the $n$-dimensional Hausdorff measure defined with respect to the metric of the ambient space $l^N_\infty$.}


%

To compare our results with this conjecture, note that we can consider a $(1+\epsilon)$-bilipschitz embedding of $M^n$ in $l^N_\infty$ as above and then apply Theorem 1.1. As a result, we will obtain an upper
bound for $W^{l^\infty}_{n-1}(M^n)$.
However, on the right-hand side we will have the volume or, equivalently, the $n$-dimensional Hausdorff measure with respect to the Euclidean metric on the ambient space $l^N_\infty$ instead of the desired $l^\infty$-Hausdorff measure $H_n(M^n)$. 
Note that each Euclidean ball of radius 
$r$ is contained in a cube with side length $2r$ with the same center, that is, the $l^N_\infty$ ball of radius $r$.
Recall that the contribution of a Euclidean ball (respectively, the concentric cube) of infinitesimally small radius $r$ to the volume $vol(M^n)$ (respectively, $H_n(M^n)$) will be equal to the volume of the Euclidean ball of radius $r$. Therefore, $vol(M^n)\geq H_n(M^n)$.
Hence, Theorem 1.1. is strictly weaker
than the inequality in Problem 1, and the right-hand side 
of the inequality in Theorem 1.1 can be up to $\sqrt{{N\over n}}$
larger than the right-hand side of the inequality in Problem 1.

\subsection{Systolic radii.}
%
%

Theorems 1.1 and 1.2 can be used to obtain new information about shortest non-contractible closed curves in essential manifolds:

\begin{theorem}
    Assume that a $C^1$-smooth closed manifold $M^n\subset \mathbb{R}^N$ is essential. 
   \par\noindent
    (1) There exists a non-contractible closed curve on $M^n$ contained in a cube in $\mathbb{R}^N$ with side length ${8\pi\over \sqrt{e}}(1+o(1))\ \sqrt{n}\ vol^{\frac{1}{n}}(M^n)$. Moreover, the cube has sides parallel to the coordinate axes. 
    \par\noindent
    (2) If the codimension $N-n$ is equal to one, then a non-contractible closed curve on $M^n$ is contained in a cube with side length $2\sqrt{3}\  vol^{\frac{1}{n}}(M^n)$ and sides parallel to the coordinate axes.
\end{theorem}

%

\begin{definition*}
For a nonsimply-connected submanifold $M^n\subset \mathbb{R}^N$ 
define $sysrad_{l^N_\infty}(M^n)$ as the infimum of radii of metric balls
in the ambient $R^N$ endowed with the $l^\infty$-metric that contain a non-contractible closed curve in $M^n$. Similarly,
$sysrad_{Eucl}(M^n)$ denotes the infimum of radii of Euclidean balls in $\mathbb{R}^N$ that contain a non-contractible closed curve in $M^n$. More generally, if 
$M^n$ is a submanifold of a Riemannian manifold or a Banach space $X$, $sysrad_X(M^n)$
is the infimum
of radii of metric balls of $X$ that contain a non-contractible closed curve in $X$.
\end{definition*}

Obviously, $$sysrad_{l^N_\infty}(M^n)\leq sysrad_{Eucl}(M^n)\leq {1\over 2}sys_1(M^n).$$ Therefore, if $M^n$ is a $C^1$-smooth closed essential manifold, we immediately obtain upper bound $({n!\over 2})^{1\over n}vol^{1\over n}(M^n)\leq {n\over 2}vol^{1\over n}(M^n)$ for
$sysrad_{Eucl}(M^n)$ and $sysrad_{l^N_\infty}(M^n)$.
If $M^n$ is embedded into $L^\infty(M^n)$ by means of the Kuratowski embedding, then 
$sysrad_{L^\infty(M^n)}={1\over 2}sys_1(M^n)$.


Theorem 1.5 asserts that for each essential $M^n\subset \mathbb{R}^N$ $$sysrad_{l^N_\infty}(M^n)\leq {4\pi\over\sqrt{e}}(1+o(1))\sqrt{n}vol^{1\over n}(M^n)$$ and for each essential $M^n\subset \mathbb{R}^{n+1}$ $$sysrad_{l^{n+1}_\infty}(M^n)\leq \sqrt{3}vol^{1\over n}(M^n),$$ which is much better than the estimates that follow from the known estimates for $sys_1(M^n)$.

Theorem 1.5 immediately follows from Theorems 1.1, 1.2 and the the first part of the following theorem that will be proven
in the next section:

\begin{theorem}
Assume that a manifold $M^n\subset \mathbb{R}^N$ is essential. Then 
\par\noindent
(1) $sysrad_{l^N_\infty}(M^n)\leq W^{l^\infty}_{n-1}(M^n)$;
\par\noindent
(2) $sysrad_{Eucl}(M^n)\leq 
W^{Eucl}_{n-1}(M^n)$.
\par\noindent
More generally, if an essential manifold $M^n$ is embedded in Riemannian manifold or a Banach space $X$, then $sysrad_X(M^n)\leq W_{n-1}^X(M^n).$
\end{theorem}

%

\subsection{Theorem 1.5(2) for tori and distortion of embeddings of flat tori to Euclidean spaces.} Observe that when the codimension is equal to $1$ the upper bound for the minimal side length of the cube that contains a non-contractible closed curve on, say,
a flat $n$-dimensional torus isometrically $C^1$-embedded in $\mathbb{R}^{n+1}$ is better by a factor of $\sqrt{n}$ than the optimal upper bound for the length of the shortest non-contractible curve. (Note that the Nash-Kuiper theorem implies that each Riemannian manifold diffeomorphic to $T^n$ can be isometrically embedded into $\mathbb{R}^n$ by means of a $C^1$-smooth embedding.) A possible 
explanation will be that each isometric embedding of flat $T^n$ to $\mathbb{R}^{n+1}$ with $sys_1\sim\sqrt{n}vol^{\frac{1}{n}}$ has a distortion $\geq const\ \sqrt{n}$. Indeed, a theorem of Khot and Naor ([KN]) implies that ``many" flat tori cannot be bi-Lipschitz embedded in any Euclidean space or even $l^2$ with distortion $\leq const\ \sqrt{n}$. In view of this result, the dimensionless estimate in
Theorem 1.5(2) might be true not only in the case of codimension $1$ but for all codimensions.

\begin{conjecture} There exists an absolute constant $const$, such that if $M^n\subset \mathbb{R}^N$ is essential, then
$sysrad_{l^N_\infty}(M^n)\leq const\ vol^{1\over n}(M^n).$

\end{conjecture}

\subsection{Some estimates for $W^{Eucl}_{n-1}(M^n)$}

A natural question is 
whether or not the estimates in Theorems 1.1, 1.2 hold for 
$W_{n-1}^{Eucl}(M^n)$,
where the distance is calculated with respect to the (Euclidean) metric on the ambient $\mathbb{R}^N$.  As Euclidean balls are smaller than $l^\infty$-balls of the same radius,
$W_{n-1}^{l^\infty}(M^n)\leq W_{n-1}^{Eucl}(M^n)$. Now, observe that
one can repeat the proof of Theorem 1.2(2) from [N] to obtain the following linear in the $n$ estimate for $W_{n-1}^{Eucl}(M^n)$:

\begin{theorem}
$$W_{n-1}^{Eucl}(M^n)\leq n\ vol^{\frac{1}{n}}(M^n).$$
\end{theorem}

Combining this estimate with Theorem 1.6(2) we obtain that $sysrad_{Eucl}(M^n)\leq nvol^{1\over n}(M^n)$,
but we already saw that
using the estimate for $sys_1(M^n)$ from [N] and the
inequality $sysrad_{Eucl}(M^n)\leq {1\over 2}sys_1(M^n)$ we can do somewhat better:

\begin{proposition} If $M^n\subset \mathbb{R}^N$ is essential, then
$$sysrad_{Eucl}(M^n)\leq
({n!\over 2})^{1\over n}vol^{1\over n}(M^n)
\leq {n\over 2} vol^{\frac{1}{n}}(M^n).$$ 
\end{proposition}

\begin{conjecture} There exists an absolute constant $const$ such that for each closed Riemannian manifold $M^n\subset \mathbb{R}^N$
$$sysrad_{Eucl}(M^n)\leq const\ \sqrt{n}vol^{1\over n}(M^n).$$
\end{conjecture}
\subsection{Estimates of $W^{l^\infty}_{n-1}(M^n)$ in terms of the Hausdorff measure of $M^n$ with respect to the $l^{\infty}$ metric.} We observed above that for $M^n\subset l^N_\infty$  Theorem 1.1 immediately implies that
$$W_{n-1}^{l^\infty}(M^n)\leq const\ \sqrt{N} H^{\frac{1}{n}}_n(M^n),$$ where $H_n(M^n)$ denotes the Hausdorff measure with respect to the  $l^\infty$-metric. However, the argument of the proof of Theorem 1.2(2) in [N]
also yields the following inequality that does not involve $N$:

\begin{theorem}
$$W_{n-1}^{l^\infty}(M^n)\leq const\ n^{\frac{3}{2}}H_n^{\frac{1}{n}}(M^n).$$
\end{theorem}

The only difference between the proof of this theorem and the previous one is that one needs to replace the coarea formula
not available for $l^N_\infty$ by the coarea inequality (a.k.a. Eilenberg's inequality), which leads to an extra $\sqrt{n}$ factor
in comparison with the inequality in Theorem 1.8. Alternatively, Theorem 1.11 immediately follows from Theorem 1.4 in [N] and the obvious inequality
between the Hausdorff content $HC_n$ in [N] and the corresponding Hausdorff measure.

\subsection{ Ideas of the proofs of Theorems 1.1 and 1.2} The proofs of Theorems 1.1 and 1.2 are based on different ideas.
In order to prove Theorem 1.1 we modify the Federer-Fleming idea
of projecting a ``small" set $M^n$ to a random shift of the $(n-1)$-skeleton
of the cubic grid in $\mathbb{R}^N$ using an appropriately shifted copy of the $(N-n)$-skeleton of the cubic grid that does not intersect $M^n$. Instead, we used an optimal $\mathbb{Z}^N$-periodic foam
in $\mathbb{R}^N$ constructed in [KORW] and [AK]. We also prove
a version of the kinematic formula,  where an averaging with respect to all isometries of the ambient Euclidean spaces is replaced by an averaging with respect to all isometries of $l^\infty$ (Theorem 3.5). Then we use this formula to
construct ``small" $\mathbb{Z}^N$-periodic intersections of appropriately shifted copies of the foam that play the role
of the randomly shifted high-codimensional skeleta of the cubic grid in $\mathbb{R}^N$. The nerve of the covering of the complement of the intersection by connected components of the complements of the shifted
copies of the foam plays the role of the dual grid.

The main idea of the proof of Theorem 1.2 is to replace a shifted copy of the $1$-skeleton of the cubic grid (that would be used in the Federer-Fleming construction) by a non-periodic curvilinear homeomorphic set in the complement of $M^n$, ($n=N-1$). The intersection of this set with each cube consists
of $2N$ paths connecting a point in the complement of $M^n$ with $2N$
points on different faces of the cube. The relative isoperimetric inequality for the $N$-cube plays a crucial role in the proof
of the existence of these graphs in the complement of $M^n$ in each cube.

\section{Proof of Theorem 1.6}

\subsection{Construction of inessential covers.}
We will prove Theorem 1.6 by contradiction. We are going to assume that there exists a positive $\delta_0$ with the following property. Let $\gamma\subset M^n\subset l^N_\infty$ be an arbitrary closed curve which is contained in an open ball of radius $W^{l^\infty}_{n-1}(M^n)+\delta_0$ in $l^N_\infty$. Then $\gamma$
is contractible in $M^n$.

The definition of $W^{l^\infty}_{n-1}(M^n)$ implies that for an arbitrary small positive
$\epsilon$ there exists a continuous map $\phi_\epsilon: M^n\longrightarrow K^{n-1}\subset \mathbb{R}^N$ such that for all $x\in M^n$ $\Vert x-\phi_\epsilon(x)\Vert_{l^\infty}< W^{l^\infty}_{n-1}(M^n)+\epsilon$. As $K^{n-1}$ is at most $(n-1)$-dimensional, there exists an arbitrarily fine open covering $\{V_\alpha\}$ of $K^{n-1}$ such that each point of $K^{n-1}$ is covered by at most $n$ open sets $\{V_\alpha\}$. For each $\delta>0$ choosing
a sufficiently fine covering, we can ensure that all connected components of sets $\phi_{\epsilon}^{-1}(V_\alpha)$ form
a covering of $M^n$ of multiplicity $\leq n$ by connected open sets contained in open balls in $l^N_\infty$ of radius $W^{l^\infty}_{n-1}(M^n)+\delta$.
If $\delta < \delta_0$, our assumption implies that the closed curves in all open sets of the covering of $M^n$ are contractible in $M^n$. This makes this open cover of $M^n$ of multiplicity $\leq n$ {\it inessential} in the terminology of Lemma A.1 in [ABHK]. This lemma implies that $M^n$ is not essential, and this completes the proof of Theorem 1.6(1).

The proof of Theorem 1.6(2) is almost identical but one uses the Euclidean metric on $\mathbb{R}^N$ instead of the $l^\infty$ metric. More generally, the same proof implies that $sysrad_X(M^n)\leq W^X_{n-1}(M^n)$.

\subsection{The existence of an inessential cover of multiplicity $\le n$ of $M^n$ implies that $M^n$ is not essential.}
\subsubsection{Introduction}
In fact, Lemma A.1 in [ABHK] asserts that $M^n$ is not essential if and only if $M^n$ has an open covering of multiplicity $\leq n$
such that all closed curves contained in one set of the covering are contractible. (This lemma is similar to Theorem $14'$ in [S], and its
proof is based on the same idea.)
Here we obviously need this lemma only in one direction. The proof of this lemma in [ABHK] is short and easy to read, but for the sake of completeness, we now sketch a somewhat different, more elementary, and more detailed proof of the fact that
the existence of an open cover of $M^n$ of multiplicity $\leq n$ by inessential connected sets $U_\alpha$ implies that $M^n$ is not essential. (This proof is a more detailed version of the argument at the end of section 2 of [N]):

Consider the universal covering $\tilde M^n$ of $M^n$. Each element $U_\alpha$ of the open covering is lifted to a collection of disjoint connected sets $\{\tilde U_{\alpha\ g}\}_{g\in \pi_1(M^n)}$ indexed 
by elements of $G=\pi_1(M^n)$. (The sets are disjoint because $U_\alpha$ is inessential.) Denote by $\tilde N$ the nerve of the open covering $\{\tilde U_{\alpha \ g}\}_{\alpha, g}$ of $\tilde M^n$. As the underlying open covering has multiplicity $\leq n$, the dimension of $\tilde N$ is at most $n-1$. Consider a $G$-equivariant partition
of unity subordinate to the open covering $\{\tilde U_{\alpha\ g}\}_{\alpha, g}$ and the corresponding $G$-equivariant map
$\tilde\psi:\tilde M^n\longrightarrow \tilde N$. The group $G$ acts freely and discontinuously on $\tilde N$, and $\tilde\psi$
is a lift of a map $\psi:M^n\longrightarrow \tilde N/G$ into a CW complex $\tilde N/G$ of dimension $\leq n-1$ that induces the
monomorphism of fundamental groups. Indeed, if $\gamma$ is a closed non-contractible loop in $M^n$, it lifts into an arc with
distinct endpoints in $\tilde M^n$, which is being mapped into an arc
$\tilde\gamma$ with distinct endpoints in $\tilde N$ by $\tilde\psi$. As $\tilde\gamma$ is the covering of the loop $\psi\circ\gamma$ in $\tilde N/G$, $\psi\circ\gamma$ is not contractible.

Below, we will prove first that $\tilde N$ is simply-connected, and then that 
$\psi$ induces an epimorphism of the fundamental groups, and, therefore, an isomorphism
of the fundamental groups. This immediately implies that $M^n$ is not essential.

In section 2.2.4 we sketch a much shorter proof of the assertion that $M^n$ is not essential, where we do not need to prove that $\tilde N$
is simply-connected.

\subsubsection{Simply-connectedness of $\tilde N$} Now we will prove that $\tilde N$ is simply connected. 

Our proof will be by contradiction. Assume that $\tilde N$
is not simply-connected.
Let $\gamma$ be a non-contractible closed curve in $\tilde N$. Without loss of generality, we can
assume that $\gamma$ is a simple simplicial curve in the $1$-skeleton of $\tilde N$ that is formed by $1$-dimensional simplices $[u_iu_{i+1}]$ of $\tilde N$, $i=1,\ldots, n_0$, $u_{n_0+1}=u_1$.
We observe that the composition of $\tilde\psi: \tilde M^n\longrightarrow \tilde N$ with each closed curve in $\tilde M^n$ is contractible (as $\tilde M^n$ is simply-connected).
Although $\gamma$ might not be the image of a closed curve in $\tilde M^n$, it is not difficult to construct a homotopic curve $\beta$ that is the image of a closed curve $\alpha$ in $\tilde M^n$. Once this is accomplished, we will know that $\gamma$ is contractible. This contradiction will prove the simply-connectedness of $\tilde N$.

In order to construct $\alpha$, (and $\beta=\tilde\psi\circ\alpha$), observe
that (1) each vertex $u_i$ corresponds to a connected open set $U_{\alpha_i g_i}$ in the open cover of $\tilde M^n$, and (2) if $u_i$
and $u_{i+1}$ are connected by an edge in the nerve $\tilde N$ of the covering, then $U_{\alpha_ig_i}\bigcap U_{\alpha_{i+1}g_{i+1}}\not=\emptyset$. 
Without any loss of generality we can assume
that the partition of unity was chosen so that for each function $t_{U_{\alpha g}}$
of the partition of unity the set $\{x\in U_{\alpha g}| t_{U_{\alpha g}}(x)>0\}$ is path-connected, and if $U_{\alpha_1 g_1}\cap U_{\alpha _2 g_2}\not= 0$, then there exists a point $x$ such that both $t_{U_{\alpha_1g_1}}(x)$ and $t_{U_{\alpha_2g_2}}(x)$ are positive.
Now, for each $i$, choose
$v_i\in U_{\alpha_ig_i}$, and $w_i\in U_{\alpha_ig_i}\bigcap U_{\alpha_{i+1}g_{i+1}}$ so that $t_{U_{\alpha_ig_i}}(v_i)>0$,
$t_{U_{\alpha_ig_i}}(w_i)>0$, and
$t_{U_{\alpha_{i+1}g_{i+1}}}(w_i)>0$.
Connect $v_i$ with $w_i$ by a path $\alpha_{i1}$ in $\{x|t_{U_{\alpha_ig_i}}(x)>0\}\subset U_{\alpha_ig_i}$, and then connect $w_i$ and $v_{i+1}$ by a path
$\alpha_{i2}$ in $\{x|t_{U_{\alpha_{i+1}g_{i+1}}}(x)>0\}\subset U_{\alpha_{i+1}g_{i+1}}$. Define $\alpha_i$ as the concatenation $\alpha_{i1}*\alpha_{i2}$ of $\alpha_{i1}$ and $\alpha_{i2}$, and $\alpha$ as the concatenation $\alpha_1*\alpha_2\ldots *\alpha_{n_0}$.
To see that there exists a homotopy from $\beta=\tilde\psi(\alpha)$ to $\gamma$, first observe that $\tilde\psi(v_i)$ is in the (open) star of $u_i$, and $\tilde\psi(\alpha_i)$ is in the union of open stars $St(u_i)$ and $St(u_{i+1})$ of vertices $u_i$ and $u_{i+1}$
in $\tilde N$.
For each $i$ choose a path $\tau_i$ from $u_i$ to $\tilde\psi(v_i)$ in $St(u_i)$ in $\tilde N$. Observe that the closed curve
$\tau_i*\tilde\psi(\alpha_i)*\tau_{i+1}^{-1}*[u_{i+1}u_i]$ is contractible, as it is contained in the contractible set $St(u_i)\cup St(u_{i+1})$.
(Here $\tau_{i+1}^{-1}$ denotes the path
$\tau_{i+1}$ traveled in the opposite direction.) 

In order to see
that $St(u_i)\cup St(u_{i+1})$ is contractible, we can consider the barycentric coordinates on $\tilde N$ (given by the functions of the chosen partition of unity). For each vertex $v$ of $\tilde N$ we will now denote the corresponding barycentric coordinate by $t_v$. Observe that $U=St(u_i)\cup St(u_{i+1})$ can be characterized as the set of all points $x$ of $\tilde N$ such that $t_{u_i}(x)+t_{u_{i+1}}(x)>0$. Therefore, one can define retraction $p:U\longrightarrow [u_iu_{i+1}]$ by
the formula $p(x)={t_{u_i}(x)\over t_{u_i}(x)+t_{i+1}(x)}u_i+{t_{u_{i+1}}(x)\over t_{u_i}(x)+t_{u_{i+1}
}(x)}u_{i+1}$.
Now we can homotope $U$ to $[u_iu_{i+1}]$ moving each $x$ towards $p(x)$ along straight line segments $[xp(x)]$. 

Therefore, for each $i$ there exists a path homotopy between $\tilde\psi(\alpha_i)$
and $\tau_i^{-1}*[u_iu_{i+1}]*\tau_{i+1}$. Combining these path homotopies for all $i$, and then canceling pairs $\tau_i*\tau_i^{-1}$, we obtain the desired homotopy from $\beta$ to $\gamma$.

Thus, we proved that $\tilde N$ is simply-connected. 
\subsubsection{Proof that $\psi$ induces an isomorphism of the fundamental groups.} Now we are going to prove that the monomorphism of the fundamental groups of $M^n$ and $\tilde N/G$ induced by $\psi$
is also an epimorphism (and, therefore,
an isomorphism). Once this is established, we will know that $M^n$ is not essential since ${\rm dim}\ \tilde N/G \leq n-1$.

We are going to proceed as follows. Choose a loop $\lambda$
representing a non-trivial element $g$ in $\pi_1(\tilde N/G)=G$ and lift it to
a path $\tilde\lambda$ connecting the base point $p_{\tilde N}$ and $gp_{\tilde N}$. Without
loss of generality, we can assume that $p_{\tilde N}$ is the image of the base point $p_{\tilde M^n}$ under $\tilde\psi$. Connect $p_{\tilde M^n}$ and $gp_{\tilde M^n}$ by a path $\tilde\rho$ in $\tilde M^n$. As $\tilde N$ is simply connected, $\tilde\psi\circ\tilde\rho$ is path homotopic to $\tilde\lambda$. 
Therefore, the covering map projects $\tilde\rho$ to a loop $\rho$ in $M^n$ such that $\psi(\rho)$ is homotopic to $\lambda$.


\subsubsection{A shorter proof, where we do not need to prove the simply-connectedness of $\tilde N$.} Instead of proving that $\tilde N$ is simply-connected, if $n\geq 3$ we could
just add $2$-discs to $\tilde N$ killing all generators of $\pi_1(\tilde N)$ in a $G$-equivariant way:  
Once we attach a $2$-cell along a closed curve $\gamma$, we also attach
$2$-cells along closed curves $g\gamma$ for all $g\in G$. As a result, we would obtain
a simply-connected complex $N_0$ of dimension $\leq n-1$ and a $G$-equivariant map
$\tilde M^n\longrightarrow \tilde N\longrightarrow N_0$ that must induce the isomorphism
of the fundamental groups of $M^n$ and $N_0/G$. If $n=2$ (and $n-1=1$), even if we do not know that $\tilde N$ is simply connected, we can be sure that $\tilde N/G$ is a $1$-dimensional complex, and its fundamental group is free. A monomorphism from the fundamental group of a closed surface to a free group can exist only if the surface is
simply-connected, and, therefore, not essential.

\forget
As it was observed at the end of section 2 in [N], $\tilde N$ is simply-connected, and, therefore, there was no need
to either add $2$-cells to kill $\pi_1(\tilde N)$ in the argument above,
or consider the case $n=2$ separately. In order to see that 
$\tilde N$ is simply connected, one can proceed as follows. Let $\gamma$ be a non-contractible curve in $\tilde N$. W.l.o.g. we can
assume that $\gamma$ is simple simplicial curve in the $1$-skeleton of $\tilde N$ that is formed by $1$-dimensiinal simplices $[u_iu_{i+1}]$ of $\tilde N$, $i=1,\ldots, n_0$, $u_{n_0+1}=u_1$.
We observe that the composition of of $\tilde\psi: \tilde M^n\longrightarrow \tilde N$ with each closed curve in $\tilde M^n$ is contractible (as $\tilde M^n$ is simply-connected).
While $\gamma$ might be not the image of a closed curve in $\tilde M^n$, it is not is not difficult to construct a homotopic curve $\beta$ that is the image of a close curve $\alpha$ in $\tilde M^n$.

In order to construct $\alpha$, (and $\beta=\tilde\psi\circ\alpha$), observe
that (1) each vertex $u_i$ corresponds to a connected open set $U_{\alpha_i g_i}$ in the open cover of $\tilde M^n$, and (2) if $u_i$
and $u_{i+1}$ are connected by an edge in the nerve $\tilde N$ of the covering, then $U_{\alpha_ig_i}\bigcap U_{\alpha_{i+1}g_{i+1}}\not=\emptyset$. Now, for each $i$ choose
$v_i\in U_{\alpha_ig_i}$, and $w_i\in U_{\alpha_ig_i}\bigcap U_{\alpha_{i+1}g_{i+1}}$. Connect $v_i$ with $w_i$ by a path $\alpha_{i1}$ in
$U_{\alpha_ig_i}$, and then connect $w_i$ and $v_{i+1}$ by a path
$\alpha_{i2}$ in $U_{\alpha_{i+1}g_{i+1}}$. Define $\alpha_i$ as the concatenation $\alpha_{i1}*\alpha_{i2}$ of $\alpha_{i1}$ and $\alpha_{i2}$, and $\alpha$ as the concatenation $\alpha_1*\alpha_2\ldots *\alpha_{n_0}$.
To see that there exists a homotopy from $\beta=\tilde\psi(\alpha)$ to $\gamma$, first observe that $\tilde\psi(v_i)$ is in the (open) star of $v_i$, and $\tilde\psi(\alpha_i)$ is in the (open) star of the closed simplex $[u_iu_{i+1}]$ in $\tilde N$.
For each $i$ choose a path $\tau_i$ from $u_i$ to $\tilde\psi(v_i)$ in the star of $u_i$ in $\tilde N$. Observe that the closed curve
$\tau_i*\tilde\psi(\alpha_i)*\tau_{i+1}^{-1}*[u_{i+1}u_i]$ is contractible, as it is contained in the (contractible) (open) star of 
$[u_iu_{i+1}]$ in $\tilde N$. (Here $\tau_{i+1}^{-1}$ denotes the path
$\tau_{i+1}$ travelled in the opposite direction.) Therefore, for each $i$ there exists a path homotopy between $\tilde\psi(\alpha_i)$
and $\tau_i^{-1}*[u_iu_{i+1}]*\tau_{i+1}$. Combining these path homotopies for all $i$, and then canceling pairs $\tau_i*\tau_i^{-1}$, we obtain the desired homotopy from $\beta$ to $\gamma$.

\forgotten
\forget
We would like to prove that 
$\psi$ induces an epimorphism of the fundamental groups.
This would imply that $\psi$ is
a continuous map of $M^n$ to a CW complex $\tilde N/G$ of dimension $\leq n-1$ that induces
an isomorphism of fundamental groups, and therefore, $M^n$ is not essential.

We would first like to prove that $\tilde N$ is simply-connected.
(This fact would imply that $\pi_1(\tilde N/G)=G$, which is less
than what we want to prove.)
While we do not need to prove the simply-connectedness of $\tilde N$ separately, the proof is a simplified and easier to read version of our proof of the surjectivity of the homomorphism of the fundamental groups induced by $\psi$, and so we prefer to explain it first.

Given a closed curve $\gamma$ in $\tilde N$, we would like to find a homotopic closed curve which is the image of a closed curve in $\tilde M^n$ under $\tilde\psi$. All closed curves in $\tilde N$ that are images of closed curves in (simply-connected manifold) $\tilde M^n$ are null-homotopic, as we can use the composition of $\tilde\psi$ and a homotopy contracting the closed curve in $\tilde M^n$.
So, assume that $\gamma$ is a closed simplicial curve in $\tilde N$. W.l.o.g. we can assume that $\gamma$ is a simplicial curve in the $1$-skeleton of $\tilde N$, and is simple.
We will ``lift" individual edges (=$1$-simplices) of $\gamma$ to $\tilde M^n$, 
so that the end points match. Yet the image of the lifted curve under $\tilde\psi$ in $\tilde N$ will not be the original curve $\gamma$, but a curve homotopic to it. 

We will first lift all the vertices on our curve. Each vertex $v$ in $\gamma$ that corresponds to some open set $U_{ig}$ in the covering of $\tilde M^n$ is incident to two edges. 
The other endpoints of these edges that will be denoted $u$ and $w$ correspond to two other sets of the covering that we will denote, correspondingly, $U_{jf}$ and $U_{kh}$. Here $g,f,h$ are elements of $G$.
Pairwise intersections of $U_{ig}$ and $U_{jf}$, as well as $U_{jf}$ and $U_{kh}$ are non-empty. 
Ideally, we would like to lift $v$ into a point in the intersection of all three of these sets, but this intersection might be empty. Therefore, we will consider two lifts of $v$ to $\tilde M^n$: $v'$ that lies in the intersection of $U_{ig}$ with $U_{jf}$, and $v"$ in the intersection of $U_{ig}$ and $U_{kh}$. We now connect $v'$ and $v"$ by a path $\lambda(v)$
in (the connected set) $U_{ig}$. This path can be regarded as a ``lift" of the vertex $v$.

Consider now the edge $[vw]$ in $\gamma$ connecting $v$ with $w$ that corresponds to $U_{kh}$. We already have chosen a point $w'$ in the intersection of $U_{ig}$ and $U_{kh}$.
We would prefer to connect $v"$ with $w'$ by a path in the intersection of $U_{ig}$ and $U_{kh}$, but this intersection might be disconnected. So, we just connect $v"$ and $w'$ by a path $\mu([uv])$  in $U_{ig}$. Now combine all paths $\lambda(v)$
corresponding to all vertices of $\gamma$ and $\mu([vw])$ corresponding to all edges of $\gamma$ into a closed curve $\alpha$ in $\tilde M^n$. We now want to prove that the curve $\beta=\tilde\psi\circ\alpha$ is homotopic to our original curve $\gamma$ in $\tilde N$.

We divide $\beta$ into arcs. Each arc will be either $\tilde\psi(\lambda(v))$ for a vertex $v$ of $\gamma$, or $\mu([uv]])$ for an edge $[uv]$ of $\gamma$. As $\lambda(v)$ is contained in $U_{ig}$ that corresponds to the vertex $v$, $\tilde\psi(\lambda(v))$ is contained in the star of $v$ in $\tilde N$. (Recall, that an (open) star $Star(L)$ of a simplicial set $L$ is a simplicial complex $K$ is the union of all open simplices in $K$
that have a face in $L$. Observe that the star $Star(\sigma)\subset K$ of each closed
simplex $\sigma$ in $K$ is contractible.)
Note that the endpoints of $\tilde\psi(\lambda(v))$ can be connected with $v$ by paths $\tau_1(v),\tau_2(v)$ in $Star(v)$.
As the triangle formed by paths $\tau_1(v),\tilde\psi(\lambda(v))$, and $\tau_2(v)$ is in
$Star(v)$, it is a contractible closed curve in $Star(v)\subset\tilde N$. Therefore, there exists a path homotopy
between $\tilde\psi(\lambda(v))$ and the path formed by $\tau_1(v)$ and
$\tau_2(v)$.

To construct the desired homotopy between $\beta$ and
$\gamma$ we start from these path homotopies for all vertices
$v$ of $\gamma$. Now it remains to construct
path homotopies between edges $[vw]$ and paths $p(v,w)$ between $v$ and $w$ formed by $\tau_2(v)$, $\tilde\psi(\mu([vw]))$, and $\tau_1(w)$ (in the oppisite direction)
for all edges $[vw]$ of $\gamma$. The existence of such path homotopies
would follow from the contractiblility of the closed curves
formed by $[vw]$ an $p(v,w)$ for all edges $[vw]$. But such a closed curve is contained
in $Star([vw])$, and is, therefore, contractible in $\tilde N$.

This completes the proof of simply-connectedness of $\tilde N$.
Before moving further, observe that we can somewhat simplify
the construction of the desired homotopy
between $\beta$ and $\gamma$ as follows: We choose an orientation of $\gamma$ and for each edge $[vw]$ of $\gamma$ we consider the concatenation $\lambda(v)*\mu([vw])$ of paths $\lambda(v)$ and $\mu([vw])$ as a lift
of $[vw]$ to $\tilde M^n$. Denote this lift of $[vw]$ by $q(v,w)$. This time we
subdivide $\alpha$ into larger arcs $\tilde\psi(q(v,w))$. To construct a homotopy from $\beta$ to $\gamma$, we start from gradually ``growing" paths $\tau_1(v)*\bar\tau_1(v)$ for all vertices $v$ of $\gamma$. Here $\tau_1(v)$ denotes that same path as before that connects the endpoint $\tilde\psi(v')$ of $\tilde\psi(q(v,w))$ with $v$, and $\bar\tau(v)$ denotes
the same path travelled in the opposite direction. Now we need to be able to construct path homotopies from $\bar\tau_1(v)*\tilde\psi(q(v,w))*\tau_1(w)$ to $[vw]$.
As before, we observe that the existence of such a path homotopy
would follow from the contractibility of the closed curve $\bar\tau_1(v)*\tilde\psi(q(v,w))*\tau_1(w)*\bar{[vw]}$, and this
closed curve is contractible since it is contained in $Star([vw])$.

Now we are going to prove that the homomoprphism of
the fundamental groups of $M^n$ and $\tilde N/G$ (that are both isomorphic to $G$) is not only a monomorphism, but also an epimorphism.

We will be using the commutative square of maps between $\tilde M^n$, $M^n$, $\tilde N$, $\tilde N/G$ formed by $\tilde\psi$,
$\psi$, the covering map $\tilde M^n\longrightarrow M^n$, and  the quotient map $\tilde N\longrightarrow \tilde N/G$. We assume that base points in these four spaces
are chosen so that they are mapped one into the other by the four considered maps between these spaces. We will be using the notation $p_S$ for the base point in a space $S$, where $S$ can be each of the considered four spaces.
Consider a non-trivial
based loop $\gamma_{\tilde N/G}$ in $\tilde N/G$. Its lifting $\gamma_{\tilde N}$ to $\tilde N$ will be a path connecting the base point $p_{\tilde N}$ with
$gp_{\tilde N}$ for some $g\in G$.

Now we would like to be able to say that $\gamma_{\tilde N}$
is the image of a path in $\tilde M^n$ under $\tilde \psi$,
but this is not necessarily true.
However, we can find a ``parallel" path
$\tilde\gamma$ in $\tilde N$ that connects the same end points,
and is path-homotopic to $\gamma_{\tilde N}$.

Divide $\gamma_{\tilde N}$ into consecutive elementary
subpaths $(\gamma_{\tilde N})_i$ so that for each $i$
$(\gamma_{\tilde N/})_i$ is contained in a closed simplex $\sigma_{f(i)}$ of $\tilde N$ (for some function $f(i)$).
As the inverse image of the star of each simplex of $\tilde N$ under $\tilde\psi$ is non-empty, one can find a path $\tilde\gamma$
connecting $p_{\tilde N}$ and $gp_{\tilde N}$ in $\tilde N$ such that: (a) $\tilde\gamma$ is also divided into consecutive 
elementary subpaths $\tilde\gamma_i$ so that for each $i$
$\tilde\gamma_i$ is contained in the (open) star of the closed simplex $\sigma_{f(i)}$ and, in addition, (b)
$\tilde\gamma$ is an image of a path $\gamma_{\tilde M^n}$ in $\tilde M^n$
under $\tilde\psi$. (This part of the proof is essentially
the same as the construction of the curve $\beta$ in our proof of the simply-connectedness of $\tilde N$ above. ) As the (open) stars of all closed simplices in $\tilde N$ are
simply-connected (an even contractible), $\tilde\gamma$ is path homotopic to to $\gamma_{\tilde N/G}$. (This homotopy can be constructed in the same way as the homotopy between $\beta$ and $\gamma$ in the proof of simply-connectedness of $\tilde N$ above.) 
Denote the image of $\gamma_{\tilde M^n}$ under the covering map
$\tilde M^n\longrightarrow M^n$ by $\gamma_{M^n}$. It is easy
to see that $[\psi\circ\gamma_{M^n}]=[\gamma_{\tilde N/G}]\in \pi_1(\tilde N/G)$. (The homotopy between these loops in $\tilde N/G$ can be constructed as the composition of the quotient map $\tilde N\longrightarrow \tilde N/G$ and the path homotopy between $\tilde\gamma$ and $\gamma_{\tilde N}$ in $\tilde N$.)
\forgotten

\forget
This is a modification of Gromov's proof of the relationship between the filling radius and 
the systole of a Riemannian manifold. Let $\phi:M^n\longrightarrow K^{n-1}$ be a continuous map such that for each 
$x\in M^n$ $\Vert x-\phi(x)\Vert\leq w=W^{l^\infty}_{n-1}(M^n)+\epsilon$, where $\epsilon$ can be made arbitrarily small.
Assume that each closed curve in a $l^\infty$-ball of radius $2w$ 
is contractible. (Of course, $l^\infty$-balls of radius $r$ is the same as cubes with side length $2r$ and sides parallel to the coordinate axes.)
By compactness, this means that for some positive $\epsilon_0$ each closed curve in a $l^\infty$-ball with side length $<2w+\epsilon_0$ is contractible.
Let $f:M^n\longrightarrow K(\pi_1(M^n), 1)$ be the classifying map.
Denote by $W$ the cylinder of the map $\phi$. 
We are going to extend $f$ to a map $F$ of the whole $W$ to 
$K(\pi_1(M^n),1)$, where $M^n$ is identified with $M^n\times \{0\}\subset W$. This would lead to a contradiction as $W$ is homotopy equivalent to $K^{n-1}$, yet $M^n$ is essential, and, therefore, $f$ cannot be factored through a $(n-1)$-dimensional complex.

Consider a very fine triangulation of $M^n$. Consider the map $\theta$ of $W$ to $\mathbb{R}^N$ that sends each $\{x\}\times [0,1]$ to the straight line segment that connects $x$ and $\phi(x)$. Consider a very fine triangulation of $W$ extending the chosen triangulation of $M^n$ so that the image of every simplex of $W$ under $\theta$ is very small. 

The desired map of $W$ is constructed inductively with respect to the dimension of the skeleton of the chosen triangulation of $W$. We start by mapping each vertex $v$ from
the $0$-skeleton to $f(\tau(v))$, where $\tau(v)$ denotes a vertex of $M^n$ that is the closest to $\theta(v)$.
Each $1$-dimensional simplex $v_1v_2$
of $W$ will be mapped to $f([\tau(v_1)\tau(v_2)])$, where $[\tau(v_1)\tau(v_2)]$ denotes the shortest geodesic in $M^n$ between $\tau(v_1)$
and $\tau(v_2)$, (and, as before, $\tau(v_1),\ \tau(v_2)$ denote the chosen nearest vertices 
to $\theta(v_1),\ \theta(v_2)$). Note that $\Vert \tau(v_1)-\tau(v_2)\Vert_{l^\infty}\leq
\Vert \tau(v_1)-v_1\Vert_{l^\infty}+\Vert v_1-v_2\Vert_{l^\infty}+\Vert v_2-\tau(v_2)\Vert_{l^\infty}\leq 2w+\delta$, where $\delta$ can be made arbitrarily small.

In order to extend to the $2$-skeleton consider a $2$-dimensional simplex $\sigma=v_1v_2v_3$
in $W$. We already mapped the boundary $\partial\sigma$ first to a geodesic triangle $T=\tau(v_1)\tau(v_2)\tau(v_3)$ in $M^n$, and then using $f$ to $K(\pi_1(M^n),1)$. It is sufficient to be able to contract $\partial T$ in $M^n$. In this
case we will extend the map of $\partial\sigma\longrightarrow\partial T$ to the topological disc $\sigma$, and then will take a composition of $f$ with the resulting
map $\sigma\longrightarrow M^n$. Observe, that $\partial T$ is contained in a $l^\infty$-ball with the center at any of its three vertices and radius $\leq 2w+\delta$, where $\delta$ can be chosen to be arbitrarily small.
Choosing $\delta<\epsilon_0$ we ensure that $\partial T$ is contractible in $M^n$, and we can extend the map $F$ to all $2$-simplices of $W$. But now we can inductively extend $F$ to all simplices of every dimension $k>2$, as the obstructions
will be in a homology group with coefficients in $\pi_{k-1}(K(\pi_1(M^n),1)=\{0\}$.
We obtain a contradiction with the assumption that all closed curves in $M^n$ contained in 
each $l^\infty$-ball of radius $\leq 2w$ are contractible.
\forgotten

%
%
%
\section{Proof of Theorem 1.1.}

\subsection{Kinematic inequalities for $l^N_\infty$.}
\subsubsection{H. Federer's kinematic formula}
We will start from mentioning the kinematic formula developed by W. Blaschke, L. Santalo, S. Chern, and H. Federer. More precisely we would like to quote a particular case of Federer's version of the formula (cf. [F],
Theorem 6.11, and Definition 2.8 for $\gamma$), when the curvature measures coincide with the Hausdorff measures.
Assume that $T^l$ and $M^m$ are two compact submanifolds ( or, more generally, smooth subpolyhedra) of $\mathbb{R}^N$, and $l+m=N+k$ for a non-negative $k$.
Let $G$ denote the group of isometries of $\mathbb{R}^N$ (generated by rotations and translations). Consider the measure $\mu$ on $G$ which is the product of the Lebesque measure on $\mathbb{R}^N$ (that coinsides with the subgroup of $G$ of all translations) and the Haar measure of volume $1$ on the orthogonal group $O(N)$.
According to the kinematic formula, the integral over $G$ with respect to $\mu$ of $vol_{k}(M^m\cap gT^{l})$ is equal to 
$\gamma(l,m,k,N)vol_{l}(T^{l})vol_m(M^m)$. (Here $g$ denotes a variable element of $G$.)
The constant $\gamma(l,m,k,N)=
\frac{\Gamma({l+1\over 2})\Gamma({m+1\over 2})}{\Gamma({k+1\over 2})\Gamma({N+1\over 2})}$.

Note that in an important for us particular case, when $l=N-1$, and $m$ and $N$ are large numbers,
$l(N-1,m,m-1,N)=\sqrt{m\over N}(1+o(1))$.
\subsubsection{Towards kinematic inequalities for $l^N_\infty$.}
We are not going to use Federer's kinematic formula in this paper.
For our purposes we need an analogous formula for the case where the group of isometries of $\mathbb{R}^n$ is replaced by the group of symmetries of $l^\infty_N$. This group is a semidirect product of the group of translations isomorphic to $\mathbb{R}^N$ and the finite group generated by permutations of coordinates and reflections with respect to coordinate hyperplanes (that is,
linear transformations $r_i$ that change the $i$th coordinate $x_i$ of $x$ to $-x_i$ leaving the other coordinates intact). This finite group is isomorphic to the group of signed permutations and has order $2^N N!$.
As in the particular case of the kinematic formula mentioned above,
we are going to first average with respect to the (finite) group of signed permutations
and then integrate with respect to the Lebesgue measure on the group of translations (isomorphic to $\mathbb{R}^N$).
However, in this context the kinematic formula apparently does not hold. Instead, we will settle for a kinematic inequality
asserting that this average/integral of the $k$-volume of intersections of two Riemannian polyhedra $T^l\bigcap gM^m$, $g\in Isom(l^N_\infty)$, of dimensions $l$ and $m$ such that $l+m=N+k$
does not exceed $\gamma'(l,m,k,N)vol_l(T^l)vol_m(M^m)$. This inequality will be established only in two cases that are relevant for us in the present paper, namely, (1) $l=N-1$; and (2) $l=N-m$. Our $\gamma'$
is greater than $\gamma$ in the Federer kinematic formula, but at
least in the important for us particular cases they have the same
asymptotics for large $l, m, N$.

\subsubsection{Upper bounds for averages of volumes of projections}

A signed permutation $(\pi, \epsilon)$ consists of a permutation $\pi$ and a $N$-dimensional vector of signs $\epsilon$. Its action on an $N$-dimensional vector $b$
transforms $b$ into the vector $(\pi, \epsilon)b=(\epsilon_i\ b_{\pi(i)})_{i=1}^N$.
Denote the group of signed permutations that acts on $\mathbb{R}^N$ (or $l^N_\infty$) by ${S}_N$.

In this section we are going to prove three lemmas.
In all three lemmas we have
two linear subspaces $L_1$, $L_2$ of $\mathbb{R}^N$, where $dim\ L_2\geq dim\ L_1$.

We apply all signed $(\pi,\epsilon)$ permutations to $L_1$ and look at the orthogonal projection $P(\pi,\epsilon):(\pi,\epsilon)L_1\longrightarrow L_2$. We prove upper bounds for the average ratio
$\frac{vol(P(\pi,\epsilon)\Omega)}{vol(\Omega)}$, where $\Omega$ is a domain in $L_1$, and we average over the group $S_N$. In Lemma 3.1 $dim\ L_1=dim\ L_2=1$, and our upper bound
for the average change of the length is ${1\over\sqrt{N}}$, Lemma 3.2 generalizes Lemma 3.1 for the case, when $dim\ L_1=1,\ dim\ L_2=m$, and our upper bound is $\sqrt{m\over N}$.
Lemma 3.3 generalizes Lemma 3.1 for the case when $dim\ L_1=dim\ L_2=m$, and the upper bound is $\frac{1}{\sqrt{N\choose m}}$.

\begin{lemma}
Let $a, b\in \mathbb{R}^N$ be two unit vectors. Then the average $\frac{1}{2^N\cdot N!}\Sigma_{(\pi,\epsilon)}\Sigma_{i=1}^N \vert\epsilon_ia_ib_{\pi(i)}\vert$ of $\vert a\cdot (\pi,\epsilon)b\vert$ does not exceed $\sqrt{\frac{1}{N}}$.
\end{lemma}

\begin{proof}
    $$\frac{1}{2^N\cdot N!}\Sigma_{(\pi,\epsilon)}\vert\Sigma_i a_ib_{\pi(i)}\epsilon_i\vert\leq \frac{1}{2^N\cdot N!}\sqrt{\Sigma_{(\pi,\epsilon)}1}\sqrt{\Sigma_{(\pi,\epsilon)}(\Sigma_i a_ib_{\pi(i)}\epsilon_i)^2}=$$
    $$\sqrt{\frac{1}{2^N\cdot N!}}\sqrt{\Sigma_{(\pi,\epsilon)}\Sigma_{i,j}a_ia_jb_{\pi(i)}b_{\pi(j)}\epsilon_i\epsilon_j}.$$

    We are going to analyze the last sum. First, note that for each pair $i,j$, where $i\not= j$ and each permutation $\pi$, exactly two combination of signs $\epsilon_i, \epsilon_j$ result in $\epsilon_i\epsilon_j$ equal to $1$ , and exactly two result in
    $-1$. Therefore, all such terms disappear. When $i=j$, we get the terms $a_i^2b^2_{\pi(i)}$ as $\epsilon^2_i$ is always equal to $1$. Therefore, the summation over all possible $2^N$ combinations of signs will yield $2^N$ under the radical sign that is canceled with the same term in the denominator of the factor $\frac{1}{2^N\cdot N!}$. Thus, the expression on the right-hand side is equal to $\sqrt{\frac{1}{N!}}\sqrt{\Sigma_\pi\Sigma_i a_i^2b^2_{\pi(i)}}.$
    In the last double sum each term $a_i^2b_j^2$ appears if and only if $\pi(i)=j$, that is, for $(N-1)!$ permutations $\pi$. Therefore, this expression is equal to $$\frac{1}{\sqrt{N!}}\sqrt{(N-1)!}\sqrt{\Sigma_{i,j}a_i^2b_j^2}=\frac{1}{\sqrt{N}}\sqrt{\Sigma_i a^2_i\Sigma_jb^2_j}=\frac{1}{\sqrt N}.$$
\end{proof}

This lemma is a particular case of the next lemma:

\begin{lemma}
Let $b\in\mathbb{R}^N$ be a unit vector, $L\subset \mathbb{R}^N$ an $m$-dimensional
linear subspace. Consider the average of the lengths of orthogonal projections of $(\pi,\epsilon)b$ on $L$, where we average over the group of signed permutations $(\pi,\epsilon)$. This average does not
exceed $\frac{\sqrt{m}}{\sqrt{N}}$.
\end{lemma}

\begin{proof}
   The proof generalizes the proof of the previous lemma. Choose an orthonormal basis $a^{(i)}$, $i=1,\dots, m$ of $L$. The length of the projection of a vector $v$ to $L$ is equal to $\sqrt{\Sigma_{j=1}^m (v\cdot a^{(j)})^2}$. The average of interest for us is
   $$\frac{1}{2^N\cdot N!}\Sigma_{(\pi,\epsilon)}\sqrt{\Sigma_{j=1}^m (\Sigma_{i=1}^N a^{(j)}_ib_{\pi(i)}\epsilon_i)^2}\leq \sqrt{\frac{1}{2^N\cdot N!}}\sqrt{\Sigma_{j=1}^m\Sigma_{(\pi,\epsilon)}(\Sigma_i a^{(j)}_ib_{\pi(i)}\epsilon_i)^2}.$$

   Now we perform the summations with respect to $(\pi,\epsilon)$ and with respect to $j$ under the radical sign. As before, for each $j$ the inner sum is equal to $2^N\cdot (N-1)!$, and we see that the expression on the right-hand side of the last inequality is equal to $\sqrt{\frac{m}{N}}.$
   
\end{proof}

Later in the paper, we will need the following similar lemma:

\begin{lemma} Let $L_1=span\{a_1,\ldots, a_m\}$, $L_2=span\{b_1,\ldots, b_m\}$ be two linear $m$-dimensional subspaces of $\mathbb{R}^N$, where $a_1, \ldots, a_m$ is an orthonormal basis of $L_1$, and $b_1,\ldots , b_m$ an orthonormal basis of $L_2$.
For each signed permutation $(\pi,\epsilon)$, consider the linear space $(\pi,\epsilon)L_2=
span\{ (\pi,\epsilon)b_1,\ldots , (\pi, \epsilon)b_m\}$. Denote the orthogonal projection $L_1\longrightarrow (\pi,\epsilon)L_2$ by $P(\pi, \epsilon)$. Denote the ratio of $\frac{vol(P(\pi,\epsilon)\Omega)}{vol(\Omega)}$ for any bounded domain $\Omega\subset L_1$ by $\theta(\pi,\epsilon)$. (It is easy to see that this ratio does not depend on $\Omega$.)
Then $\frac{1}{2^N\cdot N!}\Sigma_{(\pi,\epsilon)}\theta(\pi,\epsilon)\leq\frac{1}{\sqrt{N\choose m}}.$

\end{lemma}

\begin{proof}

Consider the matrix of $P(\pi,\epsilon)$ with respect to bases $\{a_i\}$ in $L_1$ and $\{(\pi,\epsilon)b_j\}$ in $(\pi,\epsilon)L_2$. Its entries are the dot products $a_i\cdot (\pi, \epsilon) b_j$. It can be expressed as the product of the $m\times N$ matrix $A$ with rows $a_i$ and $N\times m$ matrix  $(\pi,\epsilon) B$, where $B$ is the matrix with columns $b_j$, and $(\pi,\epsilon)$ acts on all columns of $B$.
The quantity of interest is $\frac{1}{2^N\cdot N!}\Sigma_{(\pi,\epsilon)}det(A\cdot (\pi, \epsilon)B)$.
As before, we apply the Cauchy-Schwarz inequality to
see that this quantity does not exceed $\sqrt{\frac{1}{2^N\cdot N!}\Sigma_{(\pi,\epsilon)}det^2(A\cdot (\pi,\epsilon)B)}$.
Now we need to verify that the average of $det^2(A\cdot (\pi,\epsilon)B)$
over the group of signed permutations $(\pi, \epsilon)$ does not exceed $\frac{1}{{N\choose m}}$.

Applying the Cauchy-Binet formula we see that $det(A\cdot (\pi, \epsilon)B))=\Sigma_S det A_S\cdot ((\pi, \epsilon)B)_S$, where $S$ runs over all ${N\choose m}$ subsets of
$[N]$ of size $m$, $A_S$ is the $m\times m$ minor of $A$ formed by its columns at indices
from $S$, and $((\pi, \epsilon)B)_S$ denotes the $m\times m$ minor of $(\pi, \epsilon)B$ formed by its {\it rows} at indices from $S$. Therefore, $$det^2(A\cdot (\pi, \epsilon)B)=
\Sigma_S det^2(A_S)det^2(((\pi, \epsilon)B)_S)+$$
$$\Sigma_{\{S, S'\vert S'\not= S\}}det (A_S)det(A_{S'})\det(((\pi,\epsilon)B)_S)det(((\pi, \epsilon)B)_{S'}).$$

Observe, that for each pair of distinct $m$-subsets $S, S'$ of $[N]=\{1,\ldots, N\}$ the average of $det(((\pi, \epsilon)B)_S)det(((\pi, \epsilon)B)_{S'})$ with respect to $(\pi, \epsilon)$ is equal to zero, and, therefore, the second term in the right hand side of the previous formula will vanish after one performs the averaging over the group of signed permutations. Indeed, find $k\in S'\setminus S$. For each permutation $\pi$ and a combination of $N-1$ signs $\epsilon_j,
 j\not= k,$ there are two choices for $\epsilon_k$, leading to two terms with the same absolute values but opposite signs that will cancel each other. Therefore,
 $\Sigma_{(\pi, \epsilon)}det(((\pi, \epsilon)B)_S)det(((\pi, \epsilon)B)_{S'})=0$.

 Now we need to estimate the average of
 $\Sigma_S det^2(A_S)det^2(((\pi, \epsilon)B)_S)$ over the group of signed permutation.
 As each term in this sum does not depend on the choice of signs $\epsilon_j$, the average is equal to the average of $\Sigma_S det^2(A_S)det^2((\pi B)_S)$ over the group of all permutations $\pi$. Observe, that for each pair of $m$-sets $S, S'\subset [N]$ the number of
 permutations $\pi$ such that $\{\pi^{-1}(i)\}_{i\in S}=S'$ is equal to $m!(N-m)!$. Therefore, the average 
 $$\frac{1}{N!}\Sigma_{\pi, S}det^2(A_S)det^2((\pi B)_S)=\frac{1}{{N\choose m}}\Sigma_{S, S'}det^2(A_S)det^2(B_{S'})=$$$$\frac{1}{{N\choose m}}\Sigma_S det(A_S(A_S)^T)\Sigma_{S'}det (B_{S'}(B_{S'})^T)=\frac{1}{{N\choose m}}det(AA^T)det(BB^T),$$
 where the last equality again follows from the Cauchy-Binet formula. Now it remains to recall that vectors $a_i$ form an orthonormal system, and, therefore, 
 $AA^T=I$, where $I$ denotes the identity matrix. Similarly, $BB^T=I$. Therefore, the average is equal to $\frac{1}{{N\choose m}}$, which completes the proof of the lemma.

\end{proof}

\subsubsection{Kinematic inequalities for $l^N_\infty$.} Finally, we state two $l^N_\infty$ kinematic inequalities. The first inequality is based on Lemma 3.2.
Consider the group $G_N$ of isometries of $l^N_\infty$ (or $\mathbb{R}^N$)
generated by the group ${S}_N$ of signed permutations and the group of translations $E_N$ 
that can be identified with $\mathbb{R}^N$. Note that the same group acts on spaces of compact submanifolds, or more generally, subpolyhedra of $\mathbb{R}^N$ of every dimension. (Here and below, we assume that all faces of considered subpolyhedra are smooth.)  If $P^k\subset\mathbb{R}^N$ is a $k$-dimensional polyhedron, 
$s=(\pi,\epsilon)$ a signed permutation, and $x\in \mathbb{R}^N$ a vector, 
then $sP^k+x$ will be the polyhedron obtained
from $P^k$ by acting on it by $s$, and then translating the result by $x$. 
Given a polyhedron $M^l$
of dimension $l\geq N-k$, we will be interested in the average $Ave_{S_N}vol_{k+l-N}(M^l\bigcap sP^k)$ over all $s\in S_N$, and the integral over $E_N=\mathbb{R}^N$ of $Ave_{S_N}vol_{k+l-N}(M^l\bigcap (sP^k+x))$ with respect to $dx$.
Here is our result:

\begin{theorem}
(1) Let $P^{N-1}\subset \mathbb{R}^N$ and $M^m\subset \mathbb{R}^N$ be two  
polyhedra. Then 
$$\int_{E_N}Ave_{S_N}vol_{m-1}(M^m \bigcap (sP^{N-1} +x))dx\leq \frac{\sqrt m}{\sqrt N}vol_m(M^m)vol_{N-1}(P^{N-1}).$$
\par\noindent
(2) Assume that $P^{N-1}$ is an unbounded $\mathbb{Z}^N$-periodic  polyhedron. Denote the
intersection of $P^{N-1}$ with some unit cube $C$ with sides parallel to the coordinate axes
by $P_0$, and the volume of $P_0$ by $v$. (It is easy to see that $v$ does not depend on the choice of the $l^N_\infty$-metric ball.) Then
$$\int_{[0,1]^N}Ave_{S_N}vol_{m-1}(M^m \bigcap (sP^{N-1}+x))dx\leq \frac{\sqrt m}{\sqrt N}vol_m(M^m)\cdot v.$$
Therefore, there exist a signed permutation $s_0$ and some $x_0\in [0,1]^N$ such that
$$vol_{m-1}(M^m\bigcap (s_0P^{N-1}+x_0))\leq \frac{\sqrt m}{\sqrt N}vol_m(M^m)\cdot v.$$
\par\noindent
(3) Assume that $M^m$ and $P^{N-1}$ are both unbounded $\mathbb{Z}^N$-periodic polyhedra. Then there
exists a signed permutation $s_0$ and $x_0\in [0,1]^N$ such that for each 
unit cube $C$ with sides parallel to the coordinate axes $$vol_{m-1}(M^m\bigcap (s_0P^{N-1}+x_0)\bigcap C)\leq \frac{\sqrt m}{\sqrt N}vol_m(M^m\bigcap C)vol_{N-1}(P^{N-1}\bigcap C).$$
\end{theorem}

\begin{proof}
\par\noindent
(1) Observe that by additivity, it is sufficient to prove this result when both $P^{N-1}$ and
$M^m$ are the interiors of compact manifolds with boundary. 
We can also partition $P^{N-1}$ into arbitrarily small elements that are arbitrarily close to a domain in a
tangent hyperplane. Up to infinitisemally small terms, we can assume that $P^{N-1}$ is
an open disc $D^{N-1}$ of a small radius $\Delta$ in a hyperplane with a unit normal vector $b$.
Similarly, we can reduce $M^m$ to an infinitisemally small flat disc $B^m$ of radius $\delta$,
and we can assume that $\Delta >> \delta$, yet $\Delta$ is still infinitesimally small. We will consider the intersections of $sD^{N-1}+x$ with $B^m$. The subsequent integration will yield the product $vol_m(M^n)vol_{N-1}(P^{N-1})$. We will try to determine the constant that will be multiplied by this product.

Fix the angle $\alpha$ between $B^m$ and the straight line $l$ perpendicular to $D^{N-1}$. Consider the integral of $vol_{m-1}(B^m\bigcap (D^{N-1}+x))$, where $x$ varies over $\mathbb{R}^N$. We will first integrate along $l$, and then
along $L^{N-1}$. Assume that the projection of the center of $B^m$ to $L^{N-1}$ is
in the disc concentric with $D^{N-1}$ with radius $\Delta-\delta$. Let $u_0$ denote a unit vector that spans $l$. When we translate $D^{N-1}$ along $l$ and intersect it with $B^m$, $B^m$ will be sliced in parallel round $(m-1)$-dimensional discs and the integral $\int_l vol_{m-1}(B^m\bigcap (D^{N-1}+tu_0))dt$ will be equal to
$vol_m(B^m)\vert\cos\alpha\vert$. (To formally prove this claim, observe that we are integrating
the $(m-1)$-dimensional volumes of the inverse images of the function $f$ on $B^m$ defined as the projection on $l$. The norm of the gradient of $f$ is equal to $\vert\cos\alpha\vert$. Now the assertion immediately follows from the coarea formula.)
If the projection of the center of $B^m$ to $D^{N-1}$ will be outside of the disc concentric to $D^{N-1}$ of radius $\Delta+\delta$, then the integral over $l$ will be equal to $0$, since $D^{N-1}+tu_0$ will never intersect $B^m$. If the projection of the center
of $B^m$ to $D^{N-1}$ will be in the annulus between balls with radii $\Delta-\delta$ and $\Delta+\delta$ with the same center as $D^{N-1}$, the integral over $l$ will be between $0$ and $vol_m(B^m)\vert\cos\alpha\vert$. Therefore, after integrating over $L^{N-1}$ we
will get $vol_m(B^m)\vert\cos\alpha\vert\cdot vol_{N-1}(D^{N-1})(1+o(1))$, when $\frac{\delta}{\Delta}\longrightarrow 0$.
Lemma 3.2 implies that for every $\alpha$ the average of $\vert\cos\alpha\vert$ over the angles corresponding to all $2^N N!$ signed permutations applied to $D^{N-1}$ does not exceed $\sqrt{\frac{m}{N}}$. Therefore, the average over all signed permutations $s$
of the integral over all $x\in\mathbb{R}^N$ of the volume of $B^m\cap (sD^{N-1}+x)$ will not exceed $$\frac{\sqrt{m}}{\sqrt{N}}vol_m(B^m)vol_{N-1}(D^{N-1})(1+o(1)).$$ Now, integration over $M^m$ and then $P^{N-1}$ yields the desired upper bound $$\sqrt{\frac{m}{N}}vol_m(M^m)vol_{N-1}(P^{N-1}).$$

\par\noindent
(2) Observe that $P^{N-1}=\bigcup_{K\in\mathbb{Z}^N}(P_0+K)$. Therefore,
$$\int_{[0,1]^N}Ave_{S_N}vol_{m-1}(M^m\bigcap (sP_{n-1}+x))dx=\int_{E_N}Ave_{S_N}vol_{m-1}(M^m\bigcap (sP_0+x))dx.$$
Now we can apply part (1) of the theorem to $M^m$ and $P_0$.
\par\noindent
(3) Apply part (2) to $M^m\bigcap C$ instead of $M^m$, $P^{N-1}$, and $P_0=P^{N-1}\bigcap C$.
\end{proof}

The second inequality is based on Lemma 3.3. For a set $S$, we denote the number of points in $S$ by $\vert S\vert$. If $S$ is an infinite set, $\vert S\vert=\infty$.

\begin{theorem}
(1) Let $P^{N-m}\subset \mathbb{R}^N$ and $M^m\subset \mathbb{R}^N$ be two polyhedra. Then 
$$\int_{E_N}Ave_{S_N}(\vert M^m \bigcap (sP^{N-m} +x)\vert)dx\leq \frac{1}{\sqrt{N\choose m}}vol_m(M^m)vol_{N-m}(P^{N-m}).$$
\par\noindent
(2) Assume that $P^{N-m}$ is unbounded $\mathbb{Z}$-periodic polyhedron. Denote the
intersection of $P^{N-m}$ with any unit cube $C$ with sides parallel to the coordinate axes by $P_0$, and the volume of $P_0$ by $v$. (It is easy to see that $v$ does not depend on the choice of the unit cube.) Then
$$\int_{[0,1]^N}Ave_{S_N}(\vert M^m \bigcap (sP^{N-m}+x)\vert)dx\leq \frac{1}{\sqrt {N\choose m}}vol_m(M^m)\cdot v.$$
Therefore, there exist a signed permutation $s_0$ and some $x_0\in [0,1]^N$ such that
$$\vert M^m\bigcap (s_0P^{N-m}+x_0)\vert \leq \frac{1}{\sqrt {N\choose m}}vol_m(M^m)\cdot v.$$
\end{theorem}

\begin{proof}
As in the proof of Theorem 3.4(1), consider an infinitesimally small element of $P^{N-m}$
that can be considered to be almost linear and chosen as a $(N-m)$-dimensional ball $\beta$ of radius $\Delta$ in a $(N-m)$-dimensional affine space $L$ tangent to $P^{N-m}$ and an infinitesimally small linear
element $B$ in $M^m$ that can be chosen as a ball of radius $\delta$. Decompose all translations of $\beta$ into translations along $L^\perp$ and translations along $L$.
A plane $L+t$ with $t\in L^\perp$ has a non-empty intersection with $B$ if and only if
$t$ is in the image $\pi(B)$ of the orthogonal projection of $B$ onto $L^\perp$. In a generic case, the intersection $(L+t)\bigcup B$ is a point if it is not empty. The intersection of $\beta$ and $B$ is nonempty if and only if the center of $\beta+t$ in $L+t$ is $\Delta$-close to the point
$B\bigcap (L+t)$. Thus, integration over the space of all $t$ leads to the contribution
$vol_m(\pi(B))vol_{N-m}(\beta)$. Averaging over all signed permutations and applying Lemma 3.3, we see that the average does not
exceed $\sqrt{\frac{1}{{N\choose m}}}vol_m(B)vol_{N-m}(\beta)$. Now, integration over $P^{N-m}$ and $M^m$
yields the desired upper bound $\sqrt{\frac{1}{{N\choose m}}}vol_m(M^m)vol_{N-m}(P^{N-m})$.

Part (2) follows from part (1) exactly as in the previous lemma.
\end{proof}

Here is an example of an application of the previous theorem:

\begin{lemma}
    Let $M^n\subset [0,1]^N\subset \mathbb{R}^N$ be a manifold such that its volume is less than $\sqrt{\frac{1}{{N\choose n}}}$. Consider the union $U$ of all ${N\choose n}$ $(N-n)$-dimensional coordinate linear spaces in $\mathbb{R}^N$. There exists $t\in [0,1]^N$
    such that $(U+t)\bigcap M^n$ is empty.
\end{lemma}

\begin{proof}
    We are going to apply the previous theorem. Observe that for each signed permutation $s$ $sU=U$. Therefore, there is no need for the averaging with respect to the group of signed permutations $S_N$ in Theorem 3.5(2) with $P^{N-n}=U$. Observe that the volume $v$ of the intersection $U\bigcap [0,1]^N$ is equal to ${N\choose n}$. We see that if the volume of $M^n$
    is strictly less than $\frac{1}{\sqrt{{N\choose n}}}$, then for some $t_0$, $\vert M^n\bigcap (U+t_0)\vert <1$, and therefore $\vert M^n\bigcap (U+t_0)\vert=0$.
\end{proof}

\subsection{Small separators and proof of Theorem 1.1.}

\subsubsection{Spherical cubes ([KORW],[AK])}
\begin{definition*}
   Let $B^N$ be a $N$-dimensional Banach space, $m<N$ an integer, and $w$ a positive real number.
   A $(N-m-1)$-dimensional complex $S^{N-m-1}\subset B^N$ is called an $m$-separator with mesh $w$,
   if there exists a $m$-dimensional complex $K^m\subset B^N$ 
   and a
   continuous map $\phi:B^N\setminus S^{N-m-1}\longrightarrow K^m$ such that $\Vert \phi(x)-x\Vert_{B^N}\leq w$ for each $x$.
\end{definition*}
 In this paper, we will be interested in the separators in $B^N=l^N_\infty$. The basic example
 of an $m$-separator in $l^N_\infty$ with mesh $w$ is the union of all $(N-m-1)$-dimensional affine planes $L^{N-m-1}+wC$ in $l^N_\infty$, where $L^{N-m-1}$ varies over the set of all ${N\choose N-m-1}$ coordinate subspaces, and $C$ can be an arbitrary integer vector. 
 One can choose $K^m$ as the dual $m$-dimensional lattice to the $m$-separator and $\phi$ sends points in each coordinate cube with side length $w$ to its $m$-skeleton. We will call these separators {\it cubic lattice separators}.

 Observe that the cubic lattice separators are $w\mathbb{Z}^N$-periodic, and the volume of the intersection of the cubic lattice $m$-separator with each closed coordinate $N$-dimensional cube with side length $w$ 
 and vertices in $w\mathbb{Z}^N$ is equal to $2^{m+1}{N\choose m+1}w^{N-m-1}$.

 Using an appropriate rescaling, we can assume that $w=1$. We will demonstrate the existence of a 
 $m$-separator
such that the volume of its intersection with each coordinate unit cube will be less than
$(2\pi)^{m+1} \sqrt{N(N-1)\ldots (N-m)}=(2\pi)^{m+1} \sqrt{\frac{N!}{(N-m-1)!}}$. When $m$ is much smaller than $N$,
this will be a much smaller volume than the volume of the intersection of the cubic lattice separators with the same unit cube.

First, consider the case $m=0$. The required construction can be found in the paper by N. Alon and B. Klartag ([AK]). See also an earlier paper [KORW] with a more complicated
construction proving the theorem below, and a recent paper [NR].

\begin{theorem} ([AK], Theorem 3) 
     For each $N$ there exists a subset
    $A^{N-1}$ of the unit cube in $\mathbb{R}^N$ such that (1) The union of all integer translates $A^{N-1}+z$ of $A$, $z\in\mathbb{Z}^N$, divides $\mathbb{R}^N$ into subsets of $l^\infty$ diameter $<1$;
    (2) The volume of $A^{N-1}$ does not exceed $2\pi\cdot \sqrt{N}$.
\end{theorem}

Here is a very brief sketch of the proof. We start with the $1$-parametric family of level sets $S_{\lambda}=\{x\in [0,1]^N\vert \Pi_{i=1}^N \sin (\pi x_i)=\lambda\}$ for all $\lambda\in [0,1]$. The appearance of the function $\Pi_{i=1}^N\sin (\pi x_i)$ here is due to the fact that it is an eigenfunction of the Laplacian $-\Delta$ on the space of functions on the unit $N$-dimensional cube that vanish on its boundary.
Observe that $S_\lambda$ is the boundary of the domain $\Omega_{\lambda}=\{x\in [0,1]^N\vert \Pi_{i=1}^N \sin(\pi x_i)\geq\lambda\}$. Then we choose
$\lambda\in (0,1]$ that minimizes the ratio $\frac{vol_{N-1}(S_{\lambda})}{vol_N(\Omega_\lambda)}$
in this $1$-parametric family and observe that Cheeger's inequality implies that the minimal value of this ratio does not exceed $2\pi\sqrt{N}$.
The shape of $S_{\lambda}$ for the optimal value of $\lambda$ will be somewhere between the shape of the ball and the shape of the cube
with sides parallel to the coordinate axes. This fact explains the term ``spherical cubes" in the title of [KORW]. 

Identify all pairs of opposite faces in $[0,1]^N$ (so that $[0,1]^N$ becomes the $N$-torus $T^N$), and assume that for each pair of vectors $v, w\in [0,1]^N$ (regarded as $T^N$), $v+w$ means the addition modulo $1$,
so that $v+w$ will always be in $[0,1]^N$ considered as $T^N.$

For each vector $v\in [0,1]^N$ denote $S_\lambda+v$ by $S^v$, and $\Omega_\lambda+v$ by $\Omega^v$
Let $S_0=S_\lambda=S^0$. For $i=1,2,\ldots $, choose a random vector $v_i\in [0,1]^N$. For each $i\ge 1$ define $\Omega_i$ as $\bigcup_{j\leq i}\Omega^{v_i}$, and $S_i$ as $S_{i-1}\bigcup S^{v_i}\setminus (S^{v_i}\bigcap \Omega_{i-1})$. In other words, one takes the unions of domains formed by random shifts of $\Omega_\lambda$, and at each step
add the part of the boundary of the last shifted domain that is not
inside the union of already constructed domains to $S_{i-1}$.
Alon and Klartag observed that with probability $1$ there exists a finite
$k$ such that $S_k$ will separate each pair of opposite faces of the
unit cube $[0,1]^N$. In particular, each connected component of the complement of $S_k$ has $l^\infty$-diameter strictly smaller than $1$. Then
they prove that for each $k$ the expectation of the $(N-1)$-dimensional volume
of $S_k$ does not exceed $2\pi\sqrt{N}$. Therefore, there exists a sequence of vectors $v_i,\ i=1,2,\ldots$ such that $vol(S_k)\leq 2\pi\sqrt{N}$, and $S_k$ separates all pairs of opposite faces of $[0,1]^N$.

By construction, $S_k$ is in $T^N$. Therefore, it can be lifted to a $\mathbb{Z}^N$-periodic subset of $\mathbb{R}^N$ that has the desired properties. This completes the proof of the theorem.

\subsubsection{Small separators in higher codimensions.}
Now we are going to apply the previous theorem to construct a family of separators with increasing codimensions. Note that the union of all integer translates of $A^{N-1}$ in Theorem 3.8 will be a $0$-separator. Denote it
by $SEP_0^{N-1}$. The dual complex, $K^0$, can be defined as a collection of points,
where we choose one point from each connected component $C_j$, $j=1,\ldots $ of $\mathbb{R}^N\setminus SEP_0^{N-1}$.

To construct an $m$-separator $SEP_m^{N-m-1}$ for each $m$, we are going to inductively
define the $\mathbb{Z}$-periodic $SEP_m^{N-m-1}$ such that the volume of its intersection with each $l^N_\infty$-metric ball of radius $1$ does not exceed $(2\pi)^{m+1}\sqrt{N(N-1)\ldots (N-m)}$ for each $m=1,\ldots, N-1$ by intersecting $SEP_{m-1}^{N-m}$ with
an appropriately chosen copy $(SEP_0^{N-1})_m=sSEP_0^{N-1}+x$, where $s$ is a signed permutation and $x\in [0,1]^N$. (Observe that both sets are $\mathbb{Z}$-periodic).

Theorem 3.4(2) implies the existence of a signed permutation $s$ and $x\in [0,1]^N$ such the $(N-m-1)$-dimensional volume
of the part of this intersection inside any $l^N_\infty$-metric ball of radius $1$ does not exceed $$\frac{\sqrt{N-m}}{\sqrt N}\cdot (2\pi)^m\sqrt{N(N-1)\ldots (N-m+1)}\cdot 2\pi\sqrt{N}=(2\pi)^{m+1}\sqrt{N(N-1)\ldots (N-m)}.$$
Denote the chosen $s$ by $s^{(m)}$, $x$ by $x^{(m)}$, and $s^{(m)}SEP_0^{N-1}+x^{(m)}$ by $SEP^{(m)}$.

Observe that $SEP_m^{N-m-1}$ will be the intersection of $m+1$ isometric copies $(SEP_0^{N-1})_i$, $i=0,1,\ldots , m$ of $SEP_0^{N-1}$. For each $i$, consider the (open) connected
components $C_j^{(i)}$ of the complement of $(SEP_0^{N-1})_i$. Although these components
are disjoint for each fixed $i$, the collection $\{C_j^{(i)}\}_{i,j}$ for all $i\leq m$ and $j$ will be an open cover of $\mathbb{R}^N\setminus SEP_m^{N-1}$ of multiplicity $m+1$.
The nerve of this cover will be an $m$-dimensional simplicial complex $K^m$.
We can realize $K^m$ in $\mathbb{R}^N$ by choosing one point $v_{ij}$ in each open set
$C_j^{(i)}$ and using simplices with vertices 
$v_{i_1 j_1},\ldots, v_{i_kj_k}$ each time 
$C_{j_1}^{(i_1)}\bigcap\ldots \bigcap C_{j_k}^{(i_k)}\not= \emptyset$. Note that in this case, numbers $i_1,\ldots , i_k$ are pairwise distinct.
Now using a subordinate partition of unity we obtain a continuous map $\pi_m: \mathbb{R}^N\setminus SEP_m^{N-m-1}\longrightarrow K^m$.
Each point $x$ is mapped to a 
convex linear combination of at most $m+1$ points $v_i\in C_j^{(i)}$ 
for different values of $i$ such that $x\in C_j^{(i)}$. As the $l^N_\infty$-diameter of 
$C_j^{(i)}$ is less than $1$, for each $x\in \mathbb{R}^N\setminus SEP_m^{N-m-1}$, $\Vert \phi_m(x)-x\Vert_{l^N_\infty}<1$.

Rescaling to an arbitrary positive mesh $w$, we obtain the following theorem establishing the existence of {\it small separators}:

\begin{theorem}
For each $N\geq 2$, $m=0,1,\ldots , N-1$ and a positive real $w$ there exists a $w\mathbb{Z}^N$-periodic set $wSEP_m^{N-m-1}\subset\mathbb{R}^N$ and an $m$-dimensional complex $K^m\in\mathbb{R}^N$ such that:
\par\noindent
(a) The volume of the intersection of $wSEP_m^{N-m-1}$ with any cube with side length $w$ and sides parallel to the coordinates axes does not exceed $(2\pi)^{m+1}\sqrt{N(N-1)\ldots (N-m)}w^{N-m-1}$;
\par\noindent
(b) There exists a continuous map $\phi_m: \mathbb{R}^N\setminus wSEP_m^{N-m-1}\longrightarrow K^m$ such that for each $x$ $\Vert \phi_m(x)-x\Vert_{l^N_\infty}< w$.
\end{theorem}

\par\noindent{\bf Remark 1.} As the sum of dimensions of $SEP_m^{N-m-1}$ and $K^m$ is strictly less than $\mathbb{R}^N$, we can perform an arbitrarily small perturbation of $K^m$ (and $\phi_m$) to ensure that $K^m\subset \mathbb{R}^N\setminus SEP_m^{N-m-1}$, if desired.
\par\noindent{\bf Remark 2.} The term {\it small separators} in the title of the paper refers to 
$w\mathbb{Z}^N$-periodic subsets $wSEP_m^{N-m-1}\subset l^N_\infty$ that satisfy the conditions of the previous theorem.
\par\noindent
\subsubsection{Proof of Theorem 1.1.}

We provide two proofs. The first proof is simpler and does not use Lemma 3.3. The second proof is shorter, and we hope that it can become a basis for a useful algorithm for dimensionality reduction.
\par\noindent
{\bf Proof 1:} Using a rescaling, we can assume that $(2\pi)^n\sqrt{n!}vol(M^n)<1$, and
we would need to prove that $W_{n-1}^{l^\infty}(M^n)< 1$. We inductively intersect $M^n$ with copies $s^{(i)}SEP_0^{N-1}+x^{(i)}$ for signed permutations $s^{(i)}$ and shifts $x^{(i)}\in [0,1]^N$ that are
chosen for each $i$ so that the volumes of intersections do not exceed the averages provided by Theorem 3.4(2). In each step, we get a factor $2\pi\cdot \sqrt{n-i+1}$ for $i=1,\ldots , n$. After $n$ steps we are going to get $M_0=M^n\bigcap (\bigcap_{i=1}^ns^{(i)}SEP_0^{N-1}+x^{(i)})$ such that $$vol_0(M_0)\leq (2\pi)^n\sqrt{n(n-1)\ldots 2}vol_n(M^n)<1.$$ Note that $vol_0(M_0)$ is the number of points in $M_0$. Therefore, if $vol_0(M_0)<1$, then $vol_0(M_0)=0$ and $M_0$ is empty.
This means that $M^n$ is contained in the complement of $\bigcap_{i=1}^n(s^{(i)}SEP_0^{N-1}+x^{(i)})$. 
Therefore, it can be covered by the complements to $s^{(i)}SEP_0^{N-1}+x^{(i)}$, which is an open covering of multiplicity $n$, and can be continuously mapped
to the $(n-1)$-dimensional nerve of this covering. As above, this nerve can be realized
as a complex in the ambient $\mathbb{R}^N$, with a vertex in each connected component
of $\mathbb{R}^{N}\setminus (s^{(i)}SEP_0^{N-1}+x^{(i)})$ (which has $l^\infty$-diameters $<1$), and straight line simplices formed by the chosen vertices. Each point $x$ is mapped to a convex linear combination of vertices that are points in open sets of diameter $<1$ containing $x$. Therefore, the distance between $x$ and its image is less than $1$.

\par\noindent
{\bf Proof 2:} Construct $SEP_{n-1}^{N-n}, K^{n-1}, \phi_{n-1}$ for $w=1$ as in Theorem 3.9.
According to Theorem 3.5(3) the intersection of $M^n$ with $sSEP_{n-1}^{N-n}+x$ for a random signed permutation $s$ and a random shift $x\in [0,1]^N$ is less than $$\frac{1}{\sqrt{{N\choose n}}}(2\pi)^n\sqrt{\frac{N!}{(N-n)!}}vol(M^n)=(2\pi)^n\cdot \sqrt{n!}\cdot vol(M^n)<1.$$ This
means that for some $s$  and $x\in [0,1]^N$ this intersection is empty. Therefore, $M^n$
is in the complement of $s\ SEP_{n-1}^{N-n}+x$, and we can use $\phi_{n-1}$ to map
$M^n$ to $s\ K^{n-1}+x$.



\section{Proof of Theorem 1.2}

%
%
%
%

\begin{lemma}
Let $C^{n+1}=[0,1]^{n+1}$ be a unit cube such that each of its faces intersects $M^n$
transversally. Let $A=M^n\cap C^{n+1}$. If $vol(A)<\frac{1}{3}$, then there exists a connected component of $C^{n+1}\setminus A$ whose intersection with the interior of each facet of the cube has a volume greater than $\frac{1}{2}$.
\end{lemma}
\begin{proof}

Denote by $A_i$ the $i$th connected component of $A$. Observe that the set $C^{n+1}\setminus A_i$ has two connected components. Denote by $V_i$ the connected component of $C^{n+1}\setminus A_i$ whose volume is at most $\frac{1}{2}$. 

Let us first prove that the complement (in $C^{n+1}$) $U$ of the closure of $\bigcup_i V_i$ is connected.
Observe that for each $i$ there is a very small open neighborhood of $A_i$ that does not intersect $\cup_{\{j| j\not= i\}}A_j$.
Take two points $a,b$ in $U$. We would like to connect them by a path in 
$C^{n+1}\setminus A$. The existence of such a path for all pairs of points $a$ and $b$
will demonstrate that $U$ is path connected.
First, connect $a$ and $b$ by a straight line path $\gamma$. Assume that this path intersects $A_{i_1}$, $A_{i_2}$, etc. Let $p_1$ and $q_1$ be the first
and the last points of intersection with $A_{i_1}$. Now connect $p_1$ and $q_1$ by a path
$\tau$ in $A_{i_1}$, and replace the arc of $\gamma$ between $p_1$ and $q_1$ by $\tau$.
Now perform a small perturbation of the resulting path that moves $\tau$ in the direction
of the connected component of $C^{n+1}\setminus A_i$ different from $V_i$. (Note that
$\gamma$ arrives to $p_1$ from this component and leaves $q_1$ towards the same component.)
The resulting path will connect $p_1$ and $q_1$ but will not intersect $A_{i_1}$, and we did not create any new intersections. Proceeding in a similar way we can eliminate
the intersections with $A_{i_2}, A_{i_3}, \ldots$ and obtain a path between $a$ and $b$ that does not intersect $A$.

Denote by $S_{i,j}$ the intersection of $V_i$ with the $j$th facet of the cube. (Recall, that a ``facet" of a cube means a codimension one face.)  With a slight abuse of notation we shall write $A_i$, $V_i$, or $S_{i,j}$ meaning the corresponding volumes. We shall prove that the complement $U$ of the closure of $\bigcup V_i$ is the required connected component of $C^{n+1}\setminus A$. To do so it is enough to prove that $\sum_i S_{i,j}<\frac{1}{2}$ for all $j$. And this, in turn, will immediately follow from the assumption $\sum_i A_i=A<\frac{1}{3}$ and the inequality $A_i \geq \frac{2}{3}S_{i,j}$ which we are going to now prove for all $j$.

Let $V_{i,j,t}$ be the volume of the section of $V_i$ by a hyperplane $x_j=t$, where $x_j$ is the coordinate axis orthogonal to the $j$th facet of the cube. There are two cases possible, either $|V_{i,j,t}-S_{i,j}| < \frac{2}{3}S_{i,j}$ for all $0\leq t\leq 1$ or not.

Consider the first case, that $|V_{i,j,t}-S_{i,j}| < \frac{2}{3}S_{i,j}$ for all $0\leq t\leq 1$. Then $$V_i=\int_0^1 V_{i,j,t} dt > \frac{1}{3}S_{i,j}.$$ From the relative isoperimetric inequality in the cube (see the appendix in [ABBF] or an earlier paper [H]), and from $V_i\leq \frac{1}{2}$ we know that $$A_i\geq 4V_i(1-V_i)\geq 2V_i.$$ Combining the two we get that $$A_i\geq 2V_i>\frac{2}{3}S_{i,j}.$$

Consider the second case, that $|V_{i,j,t}-S_{i,j}| \geq \frac{2}{3}S_{i,j}$ for some $t$. Let $Proj_j$ be the orthogonal projection to the $j$th facet of the cube. Clearly, every point of the symmetric difference of $S_{i,j}$ and $Proj_j(V_{i,j,t})$ belongs to $Proj_j(A_i)$. The former is at least $|V_{i,j,t}-S_{i,j}| \geq \frac{2}{3}S_{i,j}$ and the latter is at most $A_i$, so once again we get $$A_i \geq \frac{2}{3}S_{i,j}.$$

\end{proof}

We are now ready to prove the theorem. Using an appropriate rescaling, we can assume that $vol(M^n)=\frac{1}{3}-\epsilon$ for an arbitrarily small positive $\epsilon$.
Consider the integer cubic lattice in the ambient $\mathbb{R}^{n+1}$. After an arbitrarily small translation of $M^n$, we may assume that it intersects all faces of all integer unit cubes transversally.

A \emph{cross} in an $(n+1)$-dimensional cube is a point in the cube's interior together with $2(n+1)$ paths connecting the point to the interior of each of the cube's facets. The interiors of the paths are disjoint. The cross intersects the interior of each facet at a single point, these points are called \emph{endpoints} of the cross.

Using $vol(M^n)<\frac{1}{3}$ and the lemma above, we can choose a cross that is disjoint with $M^n$ in each cube of the lattice. Moreover, because the connected component in the statement of the lemma intersects each facet at the volume greater than $\frac{1}{2}$, we can assume that for every facet the two endpoints of two corresponding crosses (one for each cube containing the facet) coincide.

Consider the union of all these crosses, $T$. We would like to prove that there is a retraction $\phi$ from $\mathbb{R}^{n+1}\setminus T$ to the $(n-1)$-skeleton of the lattice such that $x$ and $\phi(x)$ belong to the same cube of the lattice for all $x\in\mathbb{R}^{n+1}\setminus T$. The restriction of $\phi$ to $M^n$ will be the desired map from $M^n$ to a $(n-1)$-dimensional polyhedron with $\Vert x-\phi(x)\Vert_{l^\infty}\leq 1$ for all $x$.

If $n+1=2$, then $\phi$ is just the map that maps each connected component of $\mathbb{R}^{n+1}\setminus T$ to the single vertex of the lattice in this component.

If $n+1\geq 4$, then a general position argument implies that there exists an isotopy of $\mathbb{R}^{n+1}$ that preserves all faces of all dimensions of the lattice and which moves $T$ to the $1$-skeleton of the dual lattice. And for $T$ being the $1$-skeleton of the dual lattice, the retraction $\phi$ is obvious.

If $n+1= 3$, then the arcs in the crosses can be knotted, and this might hinder the construction of an isotopy between $T$ and the $1$-skeleton of the dual integer lattice.
Therefore, we completely abandon the construction of $T$ above, and 
present a simpler argument. We will prove that if $n+1=3$ then we can choose $T$ disjoint with $M^2$ as the image of a translation of the $1$-skeleton of the dual lattice. For such $T$, the retraction $\phi$ is again obvious.

One way of doing that is to apply Lemma 3.6 to $N=3, n=2$. Lemma 3.6 immediately
implies that if $vol(M^2)<\frac{1}{\sqrt{3}}$, which is weaker than $vol(M^2)<\frac{1}{3}$, then there exists a desired $T$ that does not intersect $M^2$.

Another way is more elementary and does not use Lemma 3.6.

For $t\in[0,1]$ denote by $X_t,Y_t$ and $Z_t$ the unions of the planes $\{x=m+t:m\in{\mathbb{Z}}\}$, $\{y=m+t:m\in{\mathbb{Z}}\}$, and $\{z=m+t:m\in{\mathbb{Z}}\}$, respectively.
By the coarea formula $$\int_0^1 vol_1(X_t\cap M^2)dt\leq vol_2(M^2)< \frac{1}{3}.$$ 
(Here and below we use subscripts to indicate the dimension of the volumes that we are considering.)
Therefore, there exists $t$ such that $vol_1(X_t\cap M^2)<\frac{1}{3}$. Denote this $X_t$ simply by $X$. 

Likewise, by the coarea formulas $$\int_0^1 vol_1(Y_t\cap M^2)dt\leq vol_2(M^2)< \frac{1}{3},$$ $$\int_0^1 vol_0(Y_t\cap X\cap M^2)dt\leq vol_1(X\cap M^2)< \frac{1}{3}.$$ Therefore, there exists $t$ such that $vol_1(Y_t\cap M^2)<\frac{2}{3}$ and $vol_0(Y_t\cap X\cap M^2)<\frac{2}{3}$ at the same time. Denote this $Y_t$ by $Y$.

Finally, by the coarea formulas $$\int_0^1 vol_0(Z_t\cap X\cap M^2)dt\leq vol_1(X\cap M^2)<\frac{1}{3}$$ $$\int_0^1 vol_0(Z_t\cap Y\cap M^2)dt\leq vol_1(Y\cap M^2)<\frac{2}{3}.$$ Therefore, there exists $t$ such that $vol_0(Z_t\cap X\cap M^2)<1$ and $vol_0(Z_t\cap Y\cap M^2)<1$ at the same time. Denote this $Z_t$ by $Z$.

The $vol_0(Y\cap X\cap M^2)$, $vol_0(Z\cap X\cap M^2)$, and $vol_0(Z\cap Y\cap M^2)$ are zero-dimensional volumes, that is, just the numbers of points in the corresponding sets. So, the fact that they are strictly less than $1$ implies that the corresponding sets $Y\cap X\cap M^2$, $Z\cap X\cap M^2$, and $Z\cap Y\cap M^2$ are empty. This means that $T:=(Y\cap X)\cup(Z\cap X)\cup(Z\cap Y)$ is as required.

Note that the previous lemma immediately implies the first part
of the following result, which is stronger than Theorem 1.2. The second part will then immediately follow from Theorem 1.6 (and is stronger than Theorem 1.5(2)).

\begin{theorem}
Let $M^n$ be a $C^1$-smooth submanifold of $\mathbb{R}^{n+1}$. Consider a cubic lattice
in $\mathbb{R}^{n+1}$ with side length $l$. Assume that the volume of the intersection of $M^n$ with any closed
cube of the $(n+1)$-dimensional cubic lattice is strictly less than $\frac{1}{3}l^n$. Then there
exists a continuous map $\phi$ of $M^n$ to the $(n-1)$-skeleton of the cubic lattice such that
$x$ and $\phi(x)$ always belong to the same closed $(n+1)$-dimensional cube of the lattice. If
$M^n$ is essential, then there exists a non-contractible closed curve on $M^n$ contained 
in an open cube with side length $2l$ and sides parallel to the coordinate axes.
\end{theorem}
{\bf Acknowledgments.} Sergey Avvakumov was partially supported by NSERC Discovery grants and the grant 765/19 from the Israel Science Foundation (ISF). The research of Alexander Nabutovsky was partially supported
by his NSERC Discovery Grant and Simons Fellowship. This research was partially done during the visit of A.N. to SLMath in the fall of 2024.

\bigskip
\address{S.A.: School of Mathematical Sciences, Tel Aviv University, Tel Aviv 69978, Israel;
s.avvakumov@gmail.com}
\par\noindent
\address{A.N.: Department of Mathematics, Bahen Centre, 40 St. George st., Rm 6290, Toronto, Ontario, M5E2S4, Canada; alex@math.toronto.edu}

\end{document}